\documentclass[a4paper]{article}
\usepackage{amsmath,amssymb,latexsym}
\usepackage[final]{graphicx}
\newcommand{\res}{\mathop{\rm res}}

\newcommand{\supp}{\mathop{\rm supp}}
\newcommand{\Ai}{\mathop{\rm Ai}}
\newcommand{\Bi}{\mathop{\rm Bi}}

\newcommand{\field}[1]{\mathbb{#1}}

\newcommand{\Z}{\field{Z}}
\newcommand{\N}{\field{N}}
\newcommand{\C}{\field{C}}

\newcommand{\A}{{\mathcal A}}

\renewcommand{\Re}{\mathop{\rm Re}}
\renewcommand{\Im}{\mathop{\rm Im}}
\newcommand{\dist}{\mathop{\rm dist}}

\newtheorem{theorem}{Theorem}[section]
\newtheorem{lemma}[theorem]{Lemma}

\newtheorem{proposition}[theorem]{Proposition}

\newtheorem{Definition}[theorem]{Definition}

\newtheorem{Remark}[theorem]{Remark}
\newenvironment{remark}{\begin{Remark}\rm}{\end{Remark}}
\newtheorem{Example}[theorem]{Example}

{\rm \trivlist \item[\hskip \labelsep{\bf Proof. }]}%
{\hspace*{\fill}$\Box$\endtrivlist}
\newenvironment{varproof}%
{\rm \trivlist \item[\hskip \labelsep{\bf Proof}]}%
{\hspace*{\fill}$\Box$\endtrivlist}
\numberwithin{equation}{section}

\title{Strong asymptotics for Jacobi polynomials with
varying nonstandard parameters}%

\author{A.B.J.\ Kuijlaars\thanks{Department of Mathematics,
Katholieke Universiteit Leuven, Celestijnenlaan 200 B, 3001 Leuven, Belgium,
email: {\tt arno@wis.kuleuven.ac.be}}
 \and A.\
Mart{\'\i}nez-Finkelshtein\thanks{University of Almer{\'\i}a and Instituto
Carlos I de F{\'\i}sica Te{\'o}rica y Computacional, Granada University
(SPAIN). Corresponding author, email: {\tt andrei@ual.es}}
 }%
\date{January 5, 2004}
\begin{document}
\maketitle

\begin{abstract}
Strong asymptotics on the whole complex plane of a sequence of
monic Jacobi polynomials $P_n^{(\alpha _n, \beta _n)}$ is studied,
assuming that
$$
\lim_{n\to\infty} \frac{\alpha_n}{n}=A\,, \qquad \lim_{n\to\infty}
\frac{\beta _n}{n}=B\,,
$$
with $A$ and $B$ satisfying $ A > -1$, $ B>-1$, $A+B < -1$. The
asymptotic analysis is based on the non-Hermitian orthogonality of
these polynomials, and uses the Deift/Zhou steepest descent
analysis for matrix Riemann-Hilbert problems. As a corollary,
asymptotic zero behavior is derived. We show that in a generic
case the zeros distribute on the set of critical trajectories
$\Gamma$ of a certain quadratic differential according to the
equilibrium measure on $\Gamma$ in an external field. However,
when either $\alpha _n$, $\beta _n$ or $\alpha _n+\beta _n$ are
geometrically close to $\Z$, part of the zeros accumulate along a
different trajectory of the same quadratic differential.
\end{abstract}



\pagestyle{myheadings} \thispagestyle{plain} \markboth{KUIJLAARS
AND MART{\'I}NEZ-FINKELSHTEIN}{ASYMPTOTICS FOR JACOBI
POLYNOMIALS}

\section{Introduction}
We consider Jacobi polynomials $P_n^{(An, Bn)}$ with varying
negative parameters $An$ and $Bn$ such that
\begin{equation} \label{ABinequalities}
     -1 < A < 0, \qquad -1 < B < 0, \qquad -2 < A + B < -1.
\end{equation}
We will obtain strong asymptotics as $n \to \infty$ of
$P_n^{(An, Bn)}(z)$ uniformly for $z$ in any region of
the complex plane and uniformly for $A$ and $B$ in compact subsets of
the set of parameter values satisfying (\ref{ABinequalities}).
Since the asymptotics is uniform in $A$ and $B$, we also
find the asymptotics for general sequences of Jacobi polynomials
$P_n^{(\alpha_n, \beta_n)}$ such that the limits
\begin{equation} \label{ABlimits}
    A = \lim_{n \to \infty} \frac{\alpha_n}{n} \qquad \mbox{and}
    \qquad B = \lim_{n \to \infty} \frac{\beta_n}{n}
\end{equation}
exist, and satisfy (\ref{ABinequalities}).
From the asymptotics of the polynomials we will also be able
to describe the limiting behavior of the zeros.

\begin{figure}[htb] \label{fig:Jacobi_cases}
\centering \includegraphics[scale=0.9]{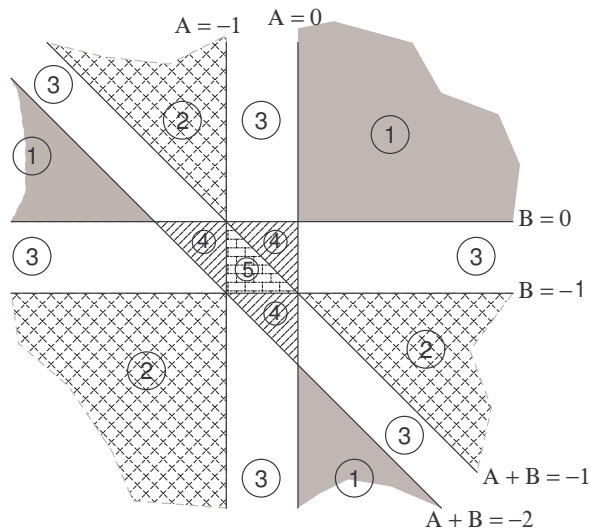}
\caption{The five different cases in the classification
of Jacobi polynomials with varying parameters
according to \cite{MR2002d:33017}.}
\end{figure}

From the point of view of behavior of zeros, the Jacobi
polynomials with varying parameters $\alpha_n$, $\beta_n$ such
that (\ref{ABlimits}) and (\ref{ABinequalities}) hold are the most
interesting general case. Indeed, Mart{\'\i}nez-Finkelshtein et al.\
\cite{MR2002d:33017} distinguish five cases depending on the
values of the limits (\ref{ABlimits}) (cf.\ Fig.\
\ref{fig:Jacobi_cases}). The first case is the case where $A, B >
0$, which corresponds to classical Jacobi polynomials with varying
positive parameters. These polynomials are orthogonal on the
interval $[-1,1]$, and as a result all their zeros are simple and
belong to the interval $(-1,1)$. The asymptotic behavior of Jacobi
polynomials with varying positive parameters is discussed in
\cite{Bosbach99,ChenIsmail,DetteStudden:95,Gawronski:91,KuijlaarsVanAssche:99,MoakSaffVarga}.
We also consider the parameter combinations $B > 0$, $A + B < -2$
and $A > 0$, $A + B < -2$ as classical. Indeed, the transformation
formula
\begin{equation} \label{transformation1}
    P_n^{(\alpha, \beta)}(x) =
    \left(\frac{1-x}{2}\right)^n P_n^{(-2n-\alpha-\beta-1, \beta)}
    \left(\frac{x+3}{x-1}\right)
\end{equation}
see \cite[\S 4 .22]{szego:1975}, expresses a Jacobi polynomial
with parameters $\alpha$ and $\beta$ satisfying
$\alpha + \beta < -2n$ and $\beta > -1$ directly in terms of a
Jacobi polynomial with positive parameters.
It follows that (\ref{transformation1}) reduces the study
of  Jacobi polynomials with varying parameters $\alpha_n$ and
$\beta_n$ such that the limits (\ref{ABlimits}) hold with $B > 0$
and $A + B < -2$ to the study of Jacobi polynomials
with varying positive parameters.
The analogous formula
\begin{equation} \label{transformation2}
    P_n^{(\alpha, \beta)}(x) =
    \left(\frac{1+x}{2}\right)^n P_n^{(\alpha, -2n-\alpha-\beta-1)}
    \left(\frac{3-x}{x+1}\right)
\end{equation}
shows similarly how to reduce the case $A > 0$
and $A+B<-2$ to the classical case.

The second case in the classification of  \cite{MR2002d:33017}
corresponds to limits $A$ and $B$ in (\ref{ABlimits}) satisfying
one of the three combinations $A < -1$, $A + B > -1$, or $B < -1$,
$A + B > -1$, or $A < -1$, $B < -1$. In this case the zeros
accumulate along an open arc in the complex plane. Their
asymptotic distribution was found in \cite{MR2002d:33017} in terms
of the equilibrium measure in an external field (cf.\
\cite{Saff:97}). The approach followed there was based on the
non-hermitian orthogonality of the Jacobi polynomials with these
parameters. See \cite{Kuijlaars03} for an overview of
non-hermitian orthogonality properties of Jacobi polynomials with
general parameters.

The remaining cases correspond to combinations of $A$ and $B$
values such that one or more of the inequalities $-1 < A < 0$,
$-1 < B< 0$, and $-2 < A+B < -1$ are satisfied.
In these cases, the zero behavior is more involved due to the possible
occurrence of multiple zeros (at $\pm 1$ only) or a zero at $\infty$
(which means a degree reduction).
To be precise, if $\alpha = -k$ is a negative integer
with $k \in \{1, \ldots, n\}$, then we have (see
\cite[formula (4.22.2)]{szego:1975}),
\begin{equation}
P_{n}^{\left( -k,\beta \right) }(z) = \frac{\Gamma(n+\beta+1)
}{\Gamma(n+\beta +1-k)}\, \frac{(n-k)!}{n!} \left(
\frac{z-1}{2}\right) ^{k}P_{n-k}^{\left( k,\beta \right) }\left(
z\right), \label{integer 1}
\end{equation}
so that $P_n^{(-k,\beta)}$ has a zero at $1$ of multiplicity $k$.
Similarly, if $\beta = -l$ with $l \in \{1, \ldots, n\}$ then
$P_n^{(\alpha, -l)}$ has a zero at $-1$ of multiplicity $l$.
A degree reduction may occur when $\alpha + \beta$ is a negative
integer, namely if $\alpha +\beta =-n -k-1$ with $k \in \{0, \ldots,n-1\}$,
then
\begin{equation}
P_{n}^{\left( \alpha ,\beta  \right) }\left( z\right)
=\frac{\Gamma(n+\alpha+1) }{\Gamma(k+\alpha+1)} \,
\frac{k!}{n!}\, P_{k}^{\left( \alpha ,\beta \right) }\left(
z\right), \label{integer 3}
\end{equation}
see \cite[Eq.\ (4.22.3)]{szego:1975}; see \S 4.22 of \cite{szego:1975}
for a more detailed discussion.
Now assume we have varying parameters $\alpha_n$, $\beta_n$ such
that the limits (\ref{ABlimits}) exist.
If $-1 < A < 0$, and if the $\alpha_n$ are integers, then
we have for each $n$ large enough, that $P_n^{(\alpha_n,\beta_n)}$
has a multiple zero at $1$. In the weak limit of the zero counting
measures this corresponds to a point mass $|A|$ at $1$.
Similarly, if $-1 < B < 0$, and if the $\beta_n$ are integers,
then we have in the limit a point mass $|B|$ at $-1$.
Finally, if $-2 < A+B< -1$ and $\alpha_n + \beta_n$ are integers,
then we have in the limit a point mass $2+A+B$ at infinity.

The classification of the remaining cases in \cite{MR2002d:33017}
depends on the number of inequalities $-1 < A<0$, $-1 < B< 0$, $-2
< A+B < -1$ that are satisfied. The third, fourth and fifth case
correspond to combination of parameters $A$ and $B$ such that
exactly one, exactly two, or exactly three, respectively, of the
inequalities are satisfied (cf.\ Fig.\ \ref{fig:Jacobi_cases}). In
these three cases the limiting behavior of zeros will be very
sensitive to the proximity of $\alpha_n$ (if $-1 < A < 0$),
$\beta_n$ (if $-1 < B < 0$) or $\alpha_n + \beta_n$ (if $-2 < A +B
< -1$) to integer values. For Laguerre polynomials the same
phenomenon was analyzed recently in
\cite{Kuijlaars/Mclaughlin:01b}.

Since all three kinds of singular behavior can occur in the fifth
case, this is the most interesting case and that is the reason why
we consider it here. The other cases can also be treated with our
methods. Fig.\ \ref{fig:Jacobi_cases} shows the behavior of zeros
which is typical for case 5. From the figure it appears that the
zeros accumulate on a contour consisting of three analytic arcs.
From our analysis below it follows that this is indeed the case,
provided that the parameters are not too close to integers.
We identify the curves as trajectories of a quadratic
differential as well as the limiting density of the zeros
on the curves, see Theorems \ref{theoremweak1} and
\ref{theoremweak2} for the exact statement.
To be able to explain the remarkable zero behavior
was the main motivation for the present work.

\begin{figure}[htb]
\centering \includegraphics[scale=0.5]{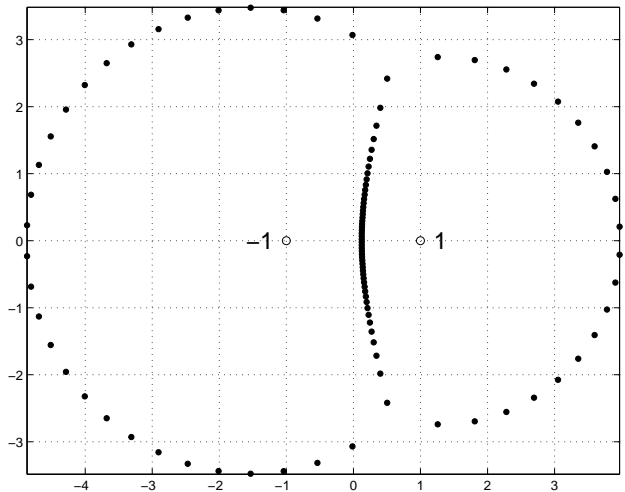} \caption{Zeros
of $P_{100}^{(\alpha ,\beta )}$ for $\alpha = -70 + 10^{-5}$,
$\beta = -80 + 10^{-5}$, which corresponds to Case 5 in Fig.\
\ref{fig:Jacobi_cases}. }\label{fig_zeros_intro}
\end{figure}

We remark that the different possibilities within
the cases 1, 2, 3, and 4 can be transformed to one
another using the transformation formulas
(\ref{transformation1}), (\ref{transformation2})
for Jacobi polynomials. It is interesting to note that
case 5 is invariant under these transformations, see \cite{MR2002d:33017}.

We also remark that the transitions between the five cases (i.e.,
$A=0$, $A=-1$, or $B=0$, $B=-1$, $A+B = -1$ or $A+B=-2$) will
present additional difficulties. These are the non-general cases,
in contrast to what we call the general cases 1--5. The zero
distribution in some of these cases has been studied by Driver,
Duren and collaborators (see also a recent survey \cite{temme03}
on the large parameter cases of the hypergeometric function). In
\cite{Driver/Duren99} the case $P_n^{(k n+1 , -n-1 )}$, $k \in
\N$, has been analyzed, corresponding to $A=k\in \N$ and $B=-1$;
this result was generalized in \cite{Duren01} using a saddle-point
method to allow $k$ to be any positive real number. Case
$P_n^{(n+b , -n-b )}$ has been studied in \cite{Driver/Moler01}.
In general, these works establish the accumulation set of the
zeros but not the limiting distribution. Trajectories of the zeros
of the Gegenbauer polynomials $P_n^{(-n-b, -n-b)}$ with fixed $n$
as $b$ varies from $-1/2$ to $-\infty$ have been described in
\cite{Driver/Duren01a}.

The rest of the paper is organized as follows. The main results
are stated in Section 2. We start defining the basic configuration
on the plane used in the description of the zero (Subsection 2.2)
and strong (Subsections 2.3--2.4) asymptotics of the polynomials.
In Section 3 we prove two technical lemmas. The cornerstone of our
approach is the matrix Riemann-Hilbert problem formulated in
Section 4; the transformations of this problem in the framework of
the Deift-Zhou steepest descent analysis (Section 5) are used in
Section 6 to prove the main results of the paper.

The Deift-Zhou steepest descent method for asymptotics of
Riemann-Hilbert problems was introduced in \cite{deift/zhou} and
applied first to orthogonal polynomials in
\cite{MR2001f:42037,MR2001g:42050}, see also \cite{MR2000g:47048}.
We use an adaptation of the method to orthogonality on curves in
the complex plane. The optimal curves are trajectories of a
quadratic differential and they were used for steepest descent
analysis of Riemann-Hilbert problems first in \cite{baik01} and
later in \cite{aptekarev02, Kamvissis03, Kuijlaars/Mclaughlin:01b,
Kuijlaars/Mclaughlin:01a, Kuijlaars03b}.

\section{Statement of results}

\subsection{Geometry of the problem}
We assume $A$ and $B$ satisfy the inequalities
(\ref{ABinequalities}) and define
\begin{equation} \label{zetapm}
    \zeta_{\pm}  =
    \frac{B^2 - A^2 \pm 4i \sqrt{(A+1)(B+1)(-A-B-1)}}
    {(A+B+2)^2}.
\end{equation}
Because of the inequalities (\ref{ABinequalities}) we have that
all factors in  the square root in (\ref{zetapm}) are positive, so
that $\zeta_+ \in \C^+=\{z\in \C:\, \Im z > 0\}$ and $\zeta_-$ is
the complex conjugate of $\zeta_+$.

Regardless of the branch of the square root and of the path of
integration we choose, the set
\begin{equation} \label{defGamma}
    \Gamma =\Gamma^{(A,B)}:= \left\{ z \in \mathbb C  :\,
        \Re \int_{\zeta_-}^z \frac{((t-\zeta_+)(t-\zeta_-))^{1/2}}{t^2-1} \, dt = 0
        \right\}
\end{equation}
is well defined, and consists of the union of the critical
trajectories of the quadratic differential (cf.\ \cite{Strebel84})
\begin{equation} \label{quaddiff}
    -\frac{(z-\zeta_-)(z-\zeta_+)}{(z^2-1)^2} \, dz^2.
\end{equation}
\begin{lemma} \label{lemma1}
We have that
$\Gamma$ is the union of three analytic arcs, which we denote by
$\Gamma_L$, $\Gamma_C$, and $\Gamma_R$. All three arcs connect the
two points $\zeta_{\pm}$ and intersect the real line in exactly
one point, in such a way that each of the intervals
$(-\infty,-1)$, $(-1,1)$, $(1,\infty)$ is cut by one of the arcs.
\end{lemma}

\begin{figure}[htb]
\centering \includegraphics[scale=0.65]{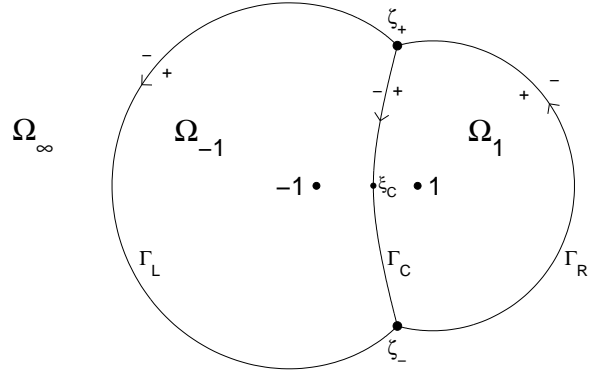}
\caption{Contour $\Gamma = \Gamma_L \cup \Gamma_C \cup \Gamma_R$
with the orientation chosen. }\label{fig_Gamma1}
\end{figure}

The contour $\Gamma$ is oriented as indicated in Fig.~\ref{fig_Gamma1}.
That is, $\Gamma_L$ and $\Gamma_C$ are oriented
from $\zeta_+$ to $\zeta_-$, and $\Gamma_R$ is oriented from
$\zeta_-$ to $\zeta_+$. The orientation of $\Gamma$ induces a $+$
and $-$ side in a neighborhood of the contour, where the $+$ side
is on the left while traversing $\Gamma$ according to its
orientation and the $-$ side is on the right. We say that a
function $f$ on $\C \setminus \Gamma$ has a boundary value
$f_+(t)$ for $t \in \Gamma \setminus \{ \zeta_+, \zeta_- \}$ if
the non-tangential limit of $f(z)$ as $z \to t$ with $z$ on the
$+$ side of $\Gamma $ exists; similarly for $f_-(t)$.

Also, we denote by $\Omega_{-1}$ and $\Omega_{1}$ the bounded
components of $\C \setminus \Gamma$ containing $-1$ and $1$,
respectively, and by $\Omega_\infty$ the unbounded component of
$\C \setminus \Gamma$ (Fig.\ \ref{fig_Gamma1}).

In what follows we write
\begin{equation} \label{defRz}
    R(z) = ((z-\zeta_+)(z-\zeta_-))^{1/2},
        \qquad z \in \mathbb C \setminus \Gamma_C,
\end{equation}
which is defined and analytic in the cut plane
$\mathbb C \setminus \Gamma_C$, such that $R(z) \sim z$
as $z \to \infty$.

We also need the critical orthogonal trajectories of the quadratic
differential (\ref{quaddiff}). These are defined by
\[ \Gamma^{\perp} = \Gamma^{\perp+} \cup \Gamma^{\perp-} \]
where
\begin{equation*} \label{gammaperp1}
\Gamma^{\perp-} = \left\{ z \in \mathbb C^-  :\,
        \Im \int_{\zeta_-}^z \frac{R(t)}{t^2-1} dt
        = 0
        \right\}
\end{equation*}
where the integration is along a path from $\zeta_-$ to $z$
in $\mathbb C^- \setminus \Gamma_C$,
and
\begin{equation*} \label{gammaperp2}
 \Gamma^{\perp+} = \left\{ z \in \mathbb C^+  :\,
        \Im \int_{\zeta_+}^z \frac{R(t)}{t^2-1} dt = 0
        \right\}
\end{equation*}
where the integration is along a path from $\zeta_+$ to $z$
in $\mathbb C^+ \setminus \Gamma_C$.

\begin{figure}[htb]
\centering \includegraphics[scale=0.65]{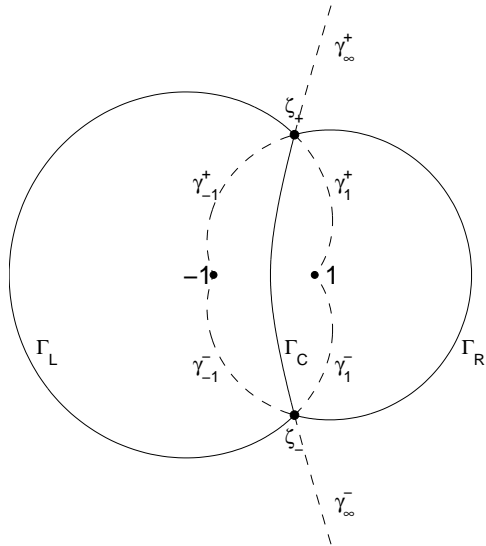}
\caption{Typical structure of the sets $\Gamma$ (solid lines) and
$\Gamma^{\perp}$ (dotted lines). }\label{fig_trajectories1}
\end{figure}

The typical structure of the orthogonal trajectories
$\Gamma^{\perp}$ is shown in Fig.~\ref{fig_trajectories1}. Three
orthogonal trajectories emanate from both $\zeta_+$ and $\zeta_-$,
ending at $1$, $-1$ and $\infty$, respectively (see the dotted
lines in Fig.\ \ref{fig_trajectories1}). We denote by
$\gamma_1^+$, $\gamma_{-1}^+$, $\gamma_{\infty}^+$ the arcs of
$\Gamma^{\perp}$ that connect $\zeta_+$ with the points $1$, $-1$,
and $\infty$, respectively; this is also the part of
$\Gamma^{\perp}$ in the upper half plane. The corresponding arcs
in the lower half plane are denoted by $\gamma_1^-$,
$\gamma_{-1}^-$ and $\gamma_{\infty}^-$, so that $\gamma_s^-$ is
the mirror image of $\gamma_s^{1}$ in the real axis, for $s \in
\{1,-1,\infty\}$.

\subsection{Weak convergence of zeros}

Then we define the absolutely continuous (a priori, complex) measure
$\mu$ on $\Gamma$ by
\begin{equation} \label{defmu}
    d\mu(z) = \frac{A + B+2}{2 \pi i} \frac{R_+(z)}{z^2-1} \,dz,
        \qquad z\in \Gamma,
\end{equation}
where $R_+$ denotes the boundary value of $R$ on the $+$-side of
$\Gamma$. (Only on $\Gamma_C$ there is a difference between the
$+$ and $-$ boundary values.) The line element $dz$ is taken
according to the orientation of $\Gamma$.

\begin{lemma} \label{lemma2}
The measure {\rm(\ref{defmu})} is positive and
\begin{equation}\label{measures_arcs1}
\mu(\Gamma_L)=1+A>0, \qquad \mu(\Gamma_C)=-1-A-B>0, \qquad
\mu(\Gamma_R)=1+B>0.
\end{equation}
In particular we have that $\mu$ is a probability measure
on $\Gamma$.
\end{lemma}

The importance of $\mu$ is shown in the following result.
\begin{theorem} \label{theoremweak1}
Let $(\alpha_n)$ and $(\beta_n)$ be two sequences such
that $\alpha_n/n \to A$ and $\beta_n/n \to B$ where $A$ and $B$
satisfy {\rm (\ref{ABinequalities})}. Suppose that
\begin{equation} \label{notclosetoZ}
    \lim_{n \to \infty} \left[ \dist(\alpha_n, \mathbb Z)\right]^{1/n} =
    \lim_{n \to \infty} \left[ \dist(\beta_n, \mathbb Z)\right]^{1/n} =
    \lim_{n \to \infty} \left[ \dist(\alpha_n+\beta_n, \mathbb Z)\right]^{1/n} = 1.
\end{equation}
Then, as $n \to \infty$, the zeros of the Jacobi polynomial
$P_n^{(\alpha_n,\beta_n)}$ accumulate on $\Gamma$ and $\mu$ is the
weak$^*$ limit of the sequence of normalized zero counting
measures.
\end{theorem}
The conditions (\ref{notclosetoZ}) imply that $\alpha_n$,
$\beta_n$, and $\alpha_n+\beta_n$ are not too close to the
integers. That such a condition is necessary is easily seen from
the case when these numbers are in fact integers (cf.\
(\ref{integer 1})).

\begin{figure}[htb]
\centering \includegraphics[scale=0.65]{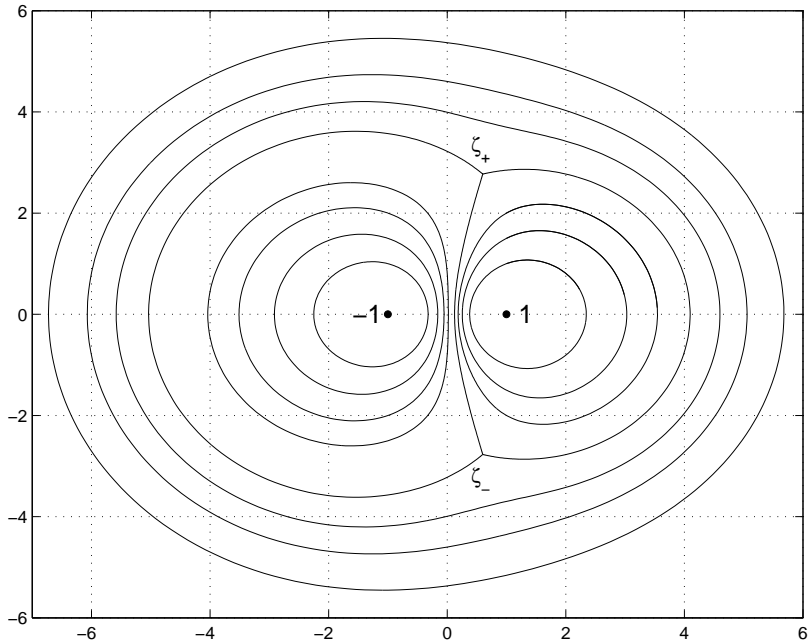}
\caption{Some trajectories of the quadratic differential
(\ref{quaddiff}), or equivalently, some level sets $\Gamma_r$, for the
values $A = -0.7$ and $B=-0.8$. } \label{trajectories}
\end{figure}

To describe the general case, we need the function
\begin{equation} \label{defphi}
    \phi(z) = \frac{A+B+2}{2} \int_{\zeta_-}^z \frac{R(t)}{t^2-1} dt,
\end{equation}
which is a multi-valued function. However, its real part is
well-defined, and we see from the definition (\ref{defGamma}) that
$\Gamma = \{ z : \Re \phi (z) = 0 \}$. For every $r$ we introduce
the level set
\begin{equation} \label{defGammar}
    \Gamma_r = \{ z \in \mathbb C :\, \Re \phi (z) = r \}.
\end{equation}
We note that by the selection of the branch in (\ref{defRz}), $\Re
\phi
> 0$ in the unbounded region $\Omega_{\infty}$ and $\Re \phi < 0$
in the two bounded regions $\Omega_{\pm 1}$. For $r > 0$, we have
that $\Gamma_r$ is a simple closed contour in $\Omega_{\infty}$,
while for $r < 0$, we have that $\Gamma_r$ consists of two simple
closed contours, one contained in $\Omega_1$ and the other in
$\Omega_{-1}$. We define for $r < 0$,
\begin{equation*} \label{defGammarmp}
    \Gamma_{r,-1} = \Gamma_r \cap \Omega_{-1},
    \qquad
    \Gamma_{r,+1} = \Gamma_r \cap \Omega_1.
\end{equation*}
We choose the positive (=counterclockwise) orientation on each of
the closed contours.
All these contours are trajectories of the quadratic differential
(\ref{quaddiff}). See Fig.~\ref{trajectories} for the trajectories.

Finally, we introduce three numbers
$r_{\alpha}$, $r_{\beta}$, and $r_{\alpha+\beta}$ and we assume that
\begin{eqnarray} \label{closetoZ1}
    \lim_{n \to \infty} \left[ \dist(\alpha_n, \mathbb Z)\right]^{1/n} & = & e^{-r_{\alpha}}, \\
     \label{closetoZ2}
    \lim_{n \to \infty} \left[ \dist(\beta_n, \mathbb Z)\right]^{1/n} & = & e^{-r_{\beta}}, \\
     \label{closetoZ3}
    \lim_{n \to \infty} \left[ \dist(\alpha_n+\beta_n, \mathbb Z)\right]^{1/n} & = & e^{-r_{\alpha+\beta}}.
\end{eqnarray}
It is easily seen that these numbers are non-negative and that the case
$r_{\alpha} = r_{\beta} = r_{\alpha+\beta} = 0$ corresponds to Theorem \ref{theoremweak1}.
It is also easily seen that at least two of the numbers $r_{\alpha}$,
$r_{\beta}$ and $r_{\alpha+\beta}$
should be equal, and if the third one is different, it will
be greater than the other two. So we distinguish four cases in the next theorem.

\begin{theorem} \label{theoremweak2}
Let $(\alpha_n)$ and $(\beta_n)$ be two sequences such
that $\alpha_n/n \to A$ and $\beta_n/n \to B$ where $A$ and $B$
satisfy {\rm (\ref{ABinequalities})}.
Suppose that there exist three numbers $r_{\alpha}$, $r_{\beta}$,
and $r_{\alpha+\beta}$ such that the
limits {\rm(\ref{closetoZ1})}, {\rm (\ref{closetoZ2})}, and {\rm (\ref{closetoZ3})}
exist. Then the following hold.
\begin{enumerate}
\item[\rm (a)] If $r_{\alpha} = r_{\beta} = r_{\alpha + \beta}$,
then the zeros of $P_n^{(\alpha_n,\beta_n)}$ accumulate on
$\Gamma$ as $n \to \infty$, and $\mu$ is the weak$^*$ limit of the
normalized zero counting measures. \item[\rm (b)] If $r_{\alpha} =
r_{\beta} < r_{\alpha+\beta}$, then the zeros of
$P_n^{(\alpha_n,\beta_n)}$ accumulate on $\Gamma_C \cup \Gamma_r$
where $r = (r_{\alpha+\beta}-r_{\alpha})/2 >0$ and
\[ \frac{A+B+2}{2\pi } \frac{R_+(z)}{z^2-1} dz, \qquad
    z \in \Gamma_C \cup \Gamma_r \]
is the weak$^*$ limit of the normalized zero counting measures.
 \item[\rm (c)] If $r_{\alpha} = r_{\alpha+\beta} < r_{\beta}$,
then the zeros of $P_n^{(\alpha_n,\beta_n)}$ accumulate on
$\Gamma_R \cup \Gamma_{r,-1}$ where $r = (r_{\alpha}-r_{\beta})/2
< 0$, and
\[ \frac{A+B+2}{2\pi i} \frac{R(z)}{z^2-1} dz, \qquad
    z \in \Gamma_R \cup \Gamma_{r,-1} \]
is the weak$^*$ limit of the normalized zero counting measures.
\item[\rm (d)] If $r_{\beta} = r_{\alpha+\beta} < r_{\alpha}$,
then the zeros of $P_n^{(\alpha_n,\beta_n)}$ accumulate
on $\Gamma_L \cup \Gamma_{r,+1}$ where $r = (r_{\beta}-r_{\alpha})/2 < 0$,
and
\[ \frac{A+B+2}{2\pi i} \frac{R(z)}{z^2-1} dz, \qquad
    z \in \Gamma_L \cup \Gamma_{r,+1} \]
is the weak$^*$ limit of the normalized zero counting measures.
\end{enumerate}
\end{theorem}
Of course the statement of Theorem \ref{theoremweak1} is a special
case of part (a) of  Theorem \ref{theoremweak2}. We choose to
mention Theorem \ref{theoremweak1} separately, since it represents
the generic case. The statements of Theorem \ref{theoremweak2} are
also valid along subsequences of $\mathbb N$, if we assume
existence of the limits (\ref{closetoZ1})--(\ref{closetoZ3}) as
$n\to \infty$ for $n$ in a subsequence $\Lambda$ of $\mathbb N$.

\begin{figure}[htb]
\centering \includegraphics[scale=0.5]{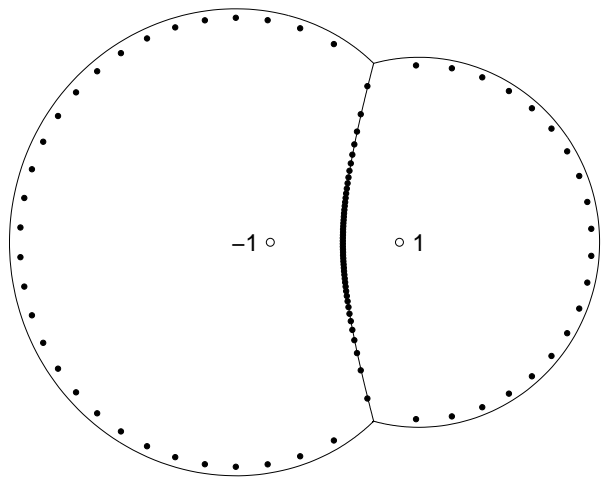} \caption{Zeros
of $P_{100}^{(\alpha ,\beta )}$ for $\alpha = -70 + 10^{-5}$,
$\beta = -80 + 10^{-5}$, together with the set $\Gamma$
corresponding to $A=-0.7$, $B=-0.8$. }\label{fig_zeros_case1a}
\end{figure}

\begin{figure}[htb]
\centering
\begin{tabular}{c@{\qquad}c}
\mbox{\includegraphics[scale=0.4]{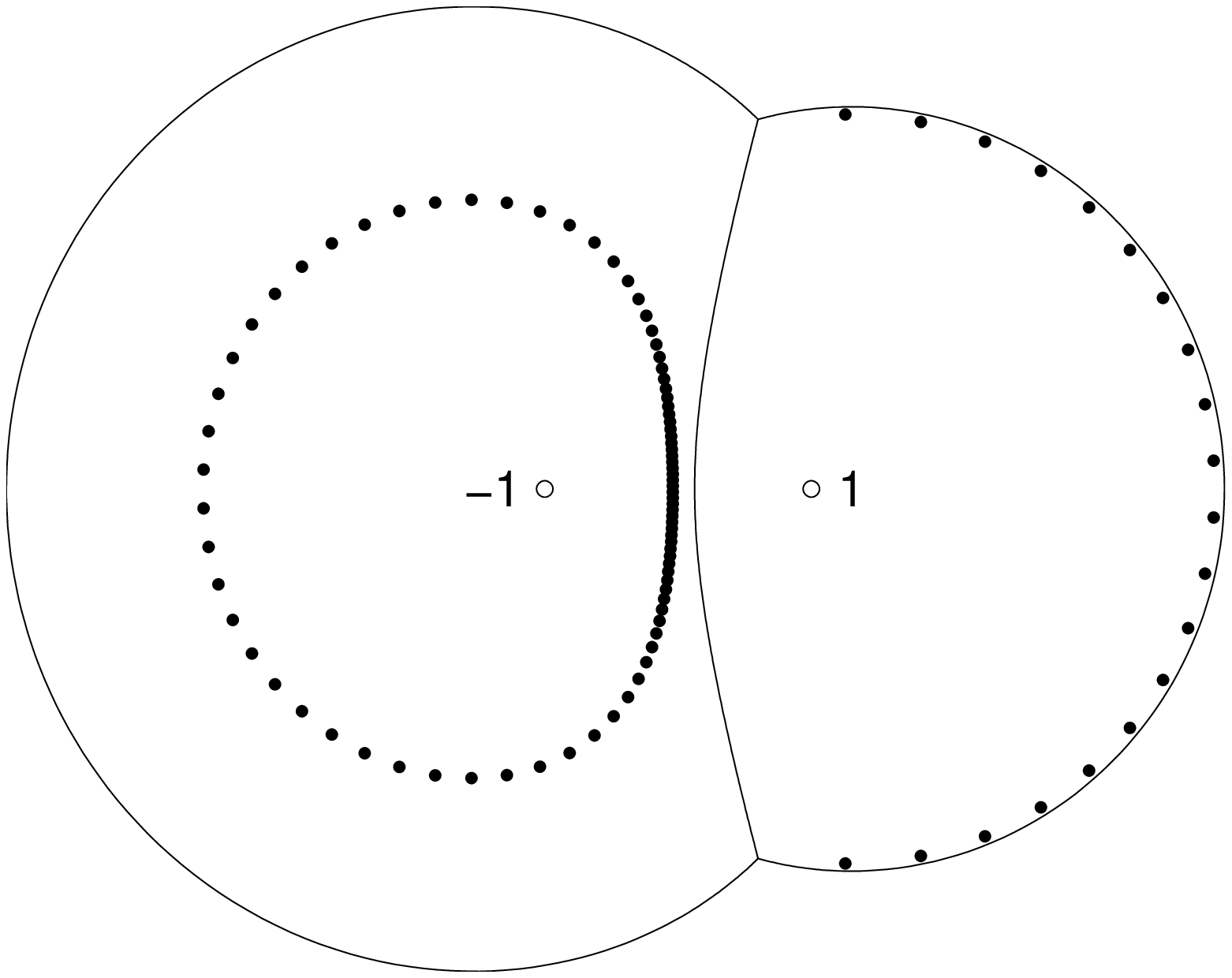}} &
\mbox{\includegraphics[scale=0.47]{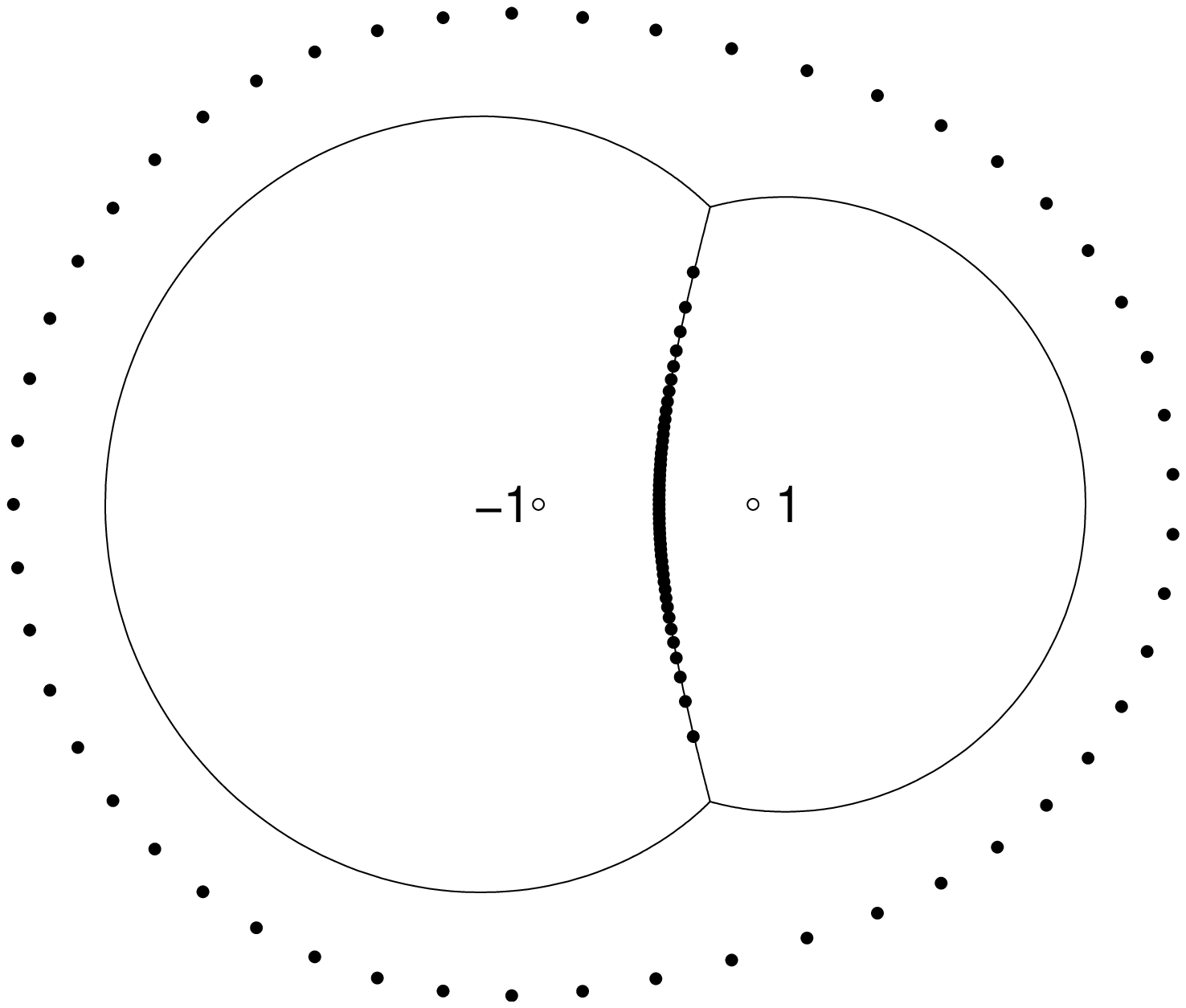}}
\end{tabular}
\caption{Zeros of $P_{100}^{(\alpha ,\beta )}$ for $\alpha = -70 +
10^{-20}$, $\beta = -80 + 10^{-30}$ (left), and $\alpha = -70 +
10^{-5}+10^{-10}$, $\beta = -80 -  10^{-5}$ (right), together with
the set $\Gamma$ corresponding to $A=-0.7$, $B=-0.8$. }
\label{fig:zerosCases23a}
\end{figure}

To illustrate the different phenomena that can happen we show some
figures (Fig.\ \ref{fig_zeros_case1a} and
Fig.\ \ref{fig:zerosCases23a}).

\begin{remark}
A general approach to the limiting zero behavior of polynomials
satisfying a non-hermitian orthogonality property has been
established in the works of Stahl \cite{Stahl:88} and
Gonchar-Rakhmanov \cite{Gonchar:87}. These authors describe the
limit distribution in terms of the equilibrium measure
in an external field on a contour satisfying a symmetry property
in $\C$. Our contour
$\Gamma$ possesses this property, but the theorems of
\cite{Stahl:88} and \cite{Gonchar:87} are not applicable: an
essential assumption in these papers is the connectedness of the complement
to the contour. Nevertheless, the measure $\mu$ from (\ref{defmu})
is the above mentioned equilibrium measure on $\Gamma$
in a certain external field.
Also the contours $\Gamma_r$ have the symmetry property
and the measures given in parts (b)--(d) of Theorem \ref{theoremweak2}
are the equilibrium measures in the external fields on these
contours. So Theorem \ref{theoremweak2} shows that in a certain sense the results
of Gonchar-Rakhmanov-Stahl are also valid for Jacobi
polynomials with varying negative parameters.
It seems likely that similar results hold in more general situations.
\end{remark}

\subsection{Strong asymptotics away from $\zeta_{\pm}$}
The weak convergence results of Theorems \ref{theoremweak1} and
\ref{theoremweak2} follow from the strong asymptotic results that we
obtain for the Jacobi polynomials. We state the result here for
the sequence $P_n^{(An, Bn)}$. We use $\widehat P_n^{(An,Bn)}$ to
denote the corresponding monic Jacobi polynomial.

\begin{figure}[htb]
\centering \includegraphics[scale=0.6]{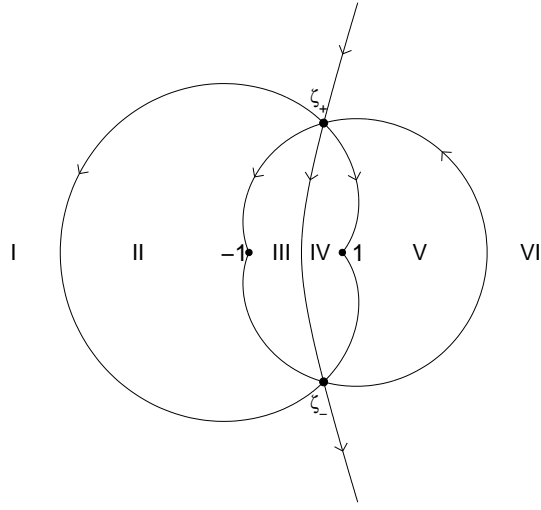}
\caption{Domains defined by trajectories $\Gamma \cup
\Gamma^\perp$. }\label{fig_regions1}
\end{figure}

Note that $\Gamma$ and $\Gamma^{\perp}$ divide the complex plane
into six domains, which we number from left to right by I, II,
III, IV, V, and VI, as shown in Fig.~\ref{fig_regions1}.

To state the asymptotic results we need to be specific about
the branches of the functions that are involved. We already
defined $\phi$ in (\ref{defphi}) as a multi-valued function.
Now we specify that
\begin{equation} \label{defphi2}
    \phi(z) = \frac{A+B+2}{2} \int_{\zeta_-}^z \frac{R(t)}{t^2-1} dt,
    \qquad z \in \mathbb C \setminus (\Gamma_C \cup \gamma_1^+ \cup
    \gamma_{-1}^+ \cup \gamma_{\infty}^+)
\end{equation}
where integration from $\zeta_-$ to $z$ is along a curve
in $\mathbb C \setminus (\Gamma_C \cup \gamma_1^+ \cup
\gamma_{-1}^+ \cup \gamma_{\infty}^+)$. Note that this
definition prevents the curve from going around the cut $\Gamma_C$
and also from going around one of the poles $\pm 1$.

Near infinity, $\phi$ behaves like
\begin{equation} \label{defc}
    \phi(z) = \frac{A+B+2}{2} \log z + c + O\left(\frac{1}{z}\right)
\end{equation}
for some constant $c$. This constant $c$ will also appear in the
asymptotic formulas below.

In our formulas we will also see fractional powers
$(z-1)^{-An/2}$ and $(z+1)^{-Bn/2}$. We will choose these
to be defined and analytic in $\mathbb C \setminus (\gamma_1^+ \cup \gamma_{\infty}^+)$
and $\mathbb C \setminus (\gamma_{-1}^+ \cup \gamma_{\infty}^+)$,
respectively, and to be positive for real $z > 1$.

Finally, we define
\begin{equation*} \label{defN11}
    N_{11}(z) =
     \frac{1}{2} \left(\left(\frac{z-\zeta_-}{z-\zeta_+}\right)^{1/4}
        + \left(\frac{z-\zeta_+}{z-\zeta_-}\right)^{1/4} \right)
\end{equation*}
and
\begin{equation*} \label{defN12}
    N_{12}(z) =
        \frac{1}{2i} \left(\left(\frac{z-\zeta_-}{z-\zeta_+}\right)^{1/4}
        - \left(\frac{z-\zeta_+}{z-\zeta_-}\right)^{1/4} \right)
\end{equation*}
which are defined and analytic in $\mathbb C \setminus \Gamma_C$.
The fourth-roots are chosen so that they approach $1$ as $z \to \infty$.
We call these functions $N_{11}$ and $N_{12}$ since they will appear
later as the corresponding entries of a matrix $N$.

Now we have all the ingredients to state our main theorem.
\begin{theorem} \label{maintheorem}
Let $A$ and $B$ satisfy {\rm (\ref{ABinequalities})}.
Then the monic Jacobi polynomials $\widehat P_n^{(An,Bn)}$ have
the following asymptotic behavior as $n \to \infty$.
\begin{enumerate}
\item[\rm (a)] For $z$ in domains I and II,
\begin{eqnarray} \nonumber
\lefteqn{\widehat P_n^{(An,Bn)}(z)
    = e^{-nc} (z-1)^{-An/2}(z+1)^{-Bn/2} } \\
    & & \label{asformleft}
        \left( e^{n\phi(z)}
        N_{11}(z)
        \left(1 + O\left(\frac{1}{n}\right)\right) \right.
        \left.
        - e^{-An\pi i} \frac{\sin (Bn\pi)}{\sin((A+B)n \pi)}
               e^{-n\phi(z)} N_{12}(z)
        \left(1 + O\left(\frac{1}{n}\right)\right) \right)
\end{eqnarray}
\item[\rm (b)] For $z$ in domain III,
\begin{equation}
\begin{split}
\widehat P_n^{(An,Bn)}(z)
    = & \ e^{-nc} (z-1)^{-An/2}(z+1)^{-Bn/2} \\
    &  \left( e^{Bn \pi i} \frac{\sin(An\pi)}{\sin((A+B)n \pi)} e^{n\phi(z)}
        N_{11}(z)
        \left(1 + O\left(\frac{1}{n}\right)\right) \right.
        \\
    & \label{asformmiddle1}
        \qquad \left. - e^{-An\pi i} \frac{\sin (Bn\pi)}{\sin((A+B)n \pi)}
               e^{-n\phi(z)} N_{12}(z)
        \left(1 + O\left(\frac{1}{n}\right)\right) \right)
\end{split}
\end{equation}
\item[\rm (c)] For $z$ in domain IV,
\begin{equation}
\begin{split}
\widehat P_n^{(An,Bn)}(z)
    = & \ e^{-nc} (z-1)^{-An/2}(z+1)^{-Bn/2}  \\
    &  \left(e^{-An\pi i} \frac{\sin (Bn\pi)}{\sin((A+B)n \pi)}
               e^{n\phi(z)} N_{11}(z)
        \left(1 + O\left(\frac{1}{n}\right)\right) \right.
        \\
    & \label{asformmiddle2}
       \qquad \left. + e^{Bn \pi i} \frac{\sin(An\pi)}{\sin((A+B)n \pi)}
        e^{-n\phi(z)} N_{12}(z)
        \left(1 + O\left(\frac{1}{n}\right)\right) \right)
\end{split}
\end{equation}
\item[\rm (d)] For $z$ in domains V and VI,
\begin{eqnarray} \nonumber
\lefteqn{\widehat P_n^{(An,Bn)}(z)
    = e^{-nc} (z-1)^{-An/2}(z+1)^{-Bn/2} } \\
    & & \label{asformright}
        \left( e^{n\phi(z)}
        N_{11}(z)
        \left(1 + O\left(\frac{1}{n}\right)\right) \right.
        \left.
        + e^{Bn \pi i} \frac{\sin(An\pi)}{\sin((A+B)n \pi)}
        e^{-n\phi(z)} N_{12}(z)
        \left(1 + O\left(\frac{1}{n}\right)\right) \right)
\end{eqnarray}
\end{enumerate}
These asymptotic formulas hold uniformly for $z$ in the indicated
domains away from the branch points, uniformly for
$A$ and $B$ in compact subsets
of the region $-1 < A < 0$, $-1 < B < 0$, $-2 < A+B < -1$,
and for values of $n$ such that $(A+B)n$ is not an integer.
\end{theorem}

\begin{remark}
One can verify that the asymptotic formulas
(\ref{asformleft})--(\ref{asformright}) agree on the boundaries
of the respective domains.
\end{remark}

\begin{remark}
The fact that the formulas (\ref{asformleft})--(\ref{asformright})
hold uniformly for $A$ and $B$ in compact
subsets of the region given by (\ref{ABinequalities}) implies that
we can allow varying values of $A$ and $B$. In particular, we can
consider two sequences $(\alpha_n)$ and $(\beta_n)$ such that the
limits (\ref{ABlimits}) exist and satisfy (\ref{ABinequalities}).
We then have asymptotic formulas for the Jacobi polynomials
$\widehat P_n^{(\alpha_n,\beta_n)}$ as in (\ref{asformleft})--(\ref{asformright})
with $A$ replaced by $A_n = \alpha_n/n$ and
$B$ replaced by $\beta_n/n$. Then we also have to realize that
$\phi$, $c$, $N_{11}$ and $N_{12}$ are going to be $n$-dependent.
Indeed, these quantities are defined using $A$ and $B$, which here
we have to replace by $A_n$ and $B_n$; we chose to state the
theorem for $\alpha _n=An$ and $\beta _n=Bn$ for the sake of
brevity of notation.
\end{remark}

\begin{remark}
The expressions between brackets in the right hand-sides of
(\ref{asformleft})--(\ref{asformright}) contain two terms that
correspond to the Liouville-Green approximation of two linearly
independent solutions of the differential equation satisfied by
the corresponding Jacobi polynomials (cf.\ \cite[Ch.\
VI]{Olver74}). In different regions of the plane and depending on
the parameters, these two terms are of comparable sizes (and then
zeros of the polynomials arise), or one of them is dominating the
other. If we assume that $An$, $Bn$, and $(A+B)n$ are not close to
integers, the expressions $\sin(An\pi)/\sin((A+B)n\pi)$ and
$\sin(Bn\pi)/\sin((A+B)n\pi)$ have moderate sizes (not too small,
not too big). In that case the dominant term is determined by $\Re
\phi$. For $z \in \Omega_{\infty}$, we have $\Re \phi(z) > 0$, and
then (\ref{asformleft}) and (\ref{asformright}) both reduce to
\begin{equation} \label{asformoutside}
    \widehat P_n^{(An,Bn)}(z)
    = e^{-nc} (z-1)^{-An/2}(z+1)^{-Bn/2}
        e^{n\phi(z)} N_{11}(z)
        \left(1 + O\left(\frac{1}{n}\right)\right)
\end{equation}
for $z$ in domains I and VI.

For $z \in \Omega_1 \cup \Omega_{-1}$ we have $\Re \phi(z) < 0$,
so that $e^{-n\phi(z)}$ dominates $e^{n\phi(z)}$ for large $n$.
Then (\ref{asformleft})--(\ref{asformright}) reduce to
\begin{equation} \label{asforminside-1}
    \widehat P_n^{(An,Bn)}(z)
    = - e^{-An\pi i-nc} \frac{\sin (Bn\pi)}{\sin((A+B)n \pi)}
        (z-1)^{-An/2}(z+1)^{-Bn/2} e^{-n\phi(z)} N_{12}(z)
        \left(1 + O\left(\frac{1}{n}\right)\right)
\end{equation}
for $z$ in domains II and III (that is, for $z \in \Omega_{-1}$), and to
\begin{equation} \label{asforminside1}
    \widehat P_n^{(An,Bn)}(z)
    = e^{Bn \pi i-nc} \frac{\sin(An\pi)}{\sin((A+B)n \pi)}
    (z-1)^{-An/2}(z+1)^{-Bn/2} e^{-n\phi(z)} N_{12}(z)
        \left(1 + O\left(\frac{1}{n}\right)\right)
\end{equation}
for $z$ in domains IV and V (that is, for $z \in \Omega_1$).

We emphasize that (\ref{asformoutside}), (\ref{asforminside-1}),
and (\ref{asforminside1}) only hold if
$An$, $Bn$, and $(A+B)n$ are not close to integers. In general
one has to use the compound asymptotic formulas
(\ref{asformleft})--(\ref{asformright}).
\end{remark}

\begin{remark}
If $An$ is an integer, then (\ref{asformmiddle2}) and
(\ref{asformright}) reduce to (\ref{asformoutside}) for $z$ in domains IV, V, and VI. Then
we see the multiple zero at $z=1$, not only because of the
factor $(z-1)^{-An/2}$, but also because
\[ \phi(z) = -\frac{A}{2} \log(z-1) + O(1) \qquad \mbox{ as } z \to 1 \]
so that
\begin{equation} \label{ephinear1}
 e^{n\phi(z)} = (z-1)^{-An/2}(1+ O(z-1)) \qquad \mbox{ as } z \to 1.
\end{equation}
So we have a zero at $z=1$ of multiplicity $-An$, as it should be.

Similar remarks apply if $Bn$ is an integer. In that case we
have a zero at $z=-1$ of multiplicity $-Bn$.
\end{remark}

\begin{remark}
If $An$ is not an integer, then $\widehat P_n^{(An,Bn)}$
does not have a zero at $z=1$. This is in agreement with formulas
(\ref{asformmiddle2}) and (\ref{asformright}) since the zero at $z=1$ due to the factor
$(z-1)^{-An/2}$ is compensated exactly by the singularity
in $e^{-n\phi(z)}$ at $z=1$, see (\ref{ephinear1}).
\end{remark}

\subsection{Strong asymptotics near $\zeta_-$.}

The asymptotic formulas (\ref{asformleft}) and (\ref{asformright})
are not valid near the branch points $\zeta_-$ and $\zeta_+$.
Near those points, there is an asymptotic formula involving
Airy functions.
We need the following particular combination of Airy functions,
depending on $A$, $B$, and $n$,
\begin{eqnarray} \label{Airycombination1} 
   \A(s;A,B,n) &= & - e^{Bn \pi i}
   \frac{\sin(An\pi)}{\sin((A+B)n\pi)}\,
        \omega \Ai(\omega s)
        +e^{-An\pi i} \frac{\sin(Bn\pi)}{\sin((A+B)n\pi)}\,
        \omega^2 \Ai(\omega^2 s) \\
        & = &
        \frac{1}{2i} \frac{\cos((A+B)n\pi)-
            \exp((B-A)n\pi i)}{\sin((A+B)n\pi)} \Ai(s)
    + \frac{1}{2i} \Bi(s),
        \label{Airycombination}
\end{eqnarray}
where $\omega = e^{2\pi i/3}$ and $\Ai$ and $\Bi$ are the usual Airy functions
\cite{Abramowitz}. Note that $\A(s;A,B,n)$ is defined for  combinations
of $A$, $B$, and $n$ that are such that $(A+B)n$ is not an integer.

\begin{theorem} \label{asymnearbranch}
Let $A$ and $B$ satisfy {\rm (\ref{ABinequalities})}. Then there
is a $\delta > 0$ such that for every $z$ with $|z-\zeta_-| <
\delta$, the monic Jacobi polynomials $\widehat P_n^{(An,Bn)}$
have the following asymptotic behavior as $n \to \infty$:
\begin{equation}
\begin{split}
     \widehat P_n^{(An, Bn)}(z)  = &
    e^{-nc} (z-1)^{-An/2} (z+1)^{-Bn/2} \sqrt{\pi} i  \\
        & \times \left[ n^{1/6}
        \left( \frac{z-\zeta_+}{z-\zeta_-} f(z)\right)^{1/4}
            \A(n^{2/3} f(z);A,B,n) \left(1+O\left(\frac{1}{n}\right)\right) \right.  \\
         &  \label{asformnearzeta-}
            \left. \qquad +
        n^{-1/6} \left(\frac{z-\zeta_+}{z-\zeta_-} f(z)\right)^{-1/4}
            \A'(n^{2/3} f(z);A,B,n) \left(1+O\left(\frac{1}{n}\right)\right) \right]
\end{split}
\end{equation}
with
\begin{equation} \label{deffz}
    f(z) = \left[ \frac{3}{2}\, \phi(z) \right]^{2/3}
\end{equation}
where the $2/3$rd root chosen is real and positive on
$\gamma_{\infty}^-$. The $O$-terms in {\rm
(\ref{asformnearzeta-})} hold uniformly for $|z-\zeta_-| < \delta$
and for $A$ and $B$ in compact subsets of the region $-1 < A < 0$,
$-1 < B < 0$, $-2 < A+B < -1$, and for values of $n$ such that
$(A+B)n$ is not an integer.
\end{theorem}

There is a similar asymptotic formula for the behavior near $\zeta_+$.

\begin{remark}
From the uniform asymptotics in Theorem \ref{asymnearbranch} it is
possible to establish a more precise behavior of the zeros of
$P_n^{(An, Bn)}$ close to the branch points $\zeta_{\pm}$. In fact,
the zeros of the function $\A$ defined in
(\ref{Airycombination1})--(\ref{Airycombination}) model the
behavior of these zeros of $P_n^{(An, Bn)}$. For instance, in the
generic case (\ref{notclosetoZ}) both terms in
(\ref{Airycombination1}) or (\ref{Airycombination}) have
approximately the same size, and $\A$ has its zeros aligned along
three curves emanating from $0$ and forming the same angle. The
situation is different in cases (b)--(d) of Theorem
\ref{theoremweak2}. For instance, if $r_{\alpha+\beta}>
r_{\alpha}= r_{\alpha}$, then the first term in
(\ref{Airycombination}) dominates the second term. But if
$r_{\alpha} \neq r_{\alpha}$ then one of the two terms in
(\ref{Airycombination1}) dominates the other. In these cases the
zeros of $\A$ behave like zeros of the dominating Airy function
and are aligned along a single curve emanating from $0$.
\end{remark}

\section{Proof of Lemmas \ref{lemma1} and \ref{lemma2}}

\begin{varproof} \textbf{of Lemma \ref{lemma1}.}
The quadratic differential (\ref{quaddiff}) has a simple zero at
$\zeta_{\pm}$ and a double pole at $\pm 1$ and at $\infty$. This
determines the local structure of the trajectories as follows, see
also \cite{baik01}, \cite[Chapter 8]{Pommerenke} or \cite[Chapter
III]{Strebel84},
\begin{enumerate}
\item[(1)] Three trajectories emanate from $\zeta_{\pm}$ at equal
angles. These are the critical trajectories.
\item[(2)] Near $\pm 1$ the trajectories are simple closed contours.
Here we use the fact that
\[ - \frac{(z-\zeta_+)(z-\zeta_-)}{(z^2-1)^2} =
    \frac{c_{\pm 1}}{(z \mp 1)^2} + O\left(\frac{1}{z \mp 1}\right) \qquad
    \textrm{as } z \to \pm 1, \]
with $c_{\pm 1} < 0$.
\item[(3)] The trajectories near $\infty$ are also simple closed contours.
This follows from the fact that in the expansion
\[ - \frac{(z-\zeta_+)(z-\zeta_-)}{(z^2-1)^2} =
    \frac{c_{\infty}}{z^2} + O\left(\frac{1}{z^3}\right)
    \qquad \textrm{as } z \to \infty, \]
we have $c_{\infty} = -1 < 0$.
\end{enumerate}

In the lower half-plane $\mathbb C^-$ there is only
the simple zero at $\zeta_-$. The other points are regular points.
This means that the three critical trajectories
that emanate from $\zeta_-$ extend to the boundary of
$\mathbb C^-$, cf.\ \cite[Lemma 8.4]{Pommerenke}.
Because of (2) and (3) and the fact that trajectories
do not intersect, the critical trajectories do not tend to infinity,
or come to $\pm 1$.
So each critical trajectories exits the lower half-plane in
a point from $\mathbb R \setminus \{-1,1\}$ and these points
are mutually distinct, say $\xi_L$, $\xi_C$, and $\xi_R$,
with $\xi_L < \xi_C < \xi_R$.
Because of the symmetry with respect to the real axis,
\[
\Re \int_{\zeta _-}^{\zeta _+} \frac{R(t)}{t^2-1}\, dt =0.
\]
Hence, the three trajectories extend into the upper half-plane as
their mirror images in $\mathbb R$, and so they continue to
$\zeta_+$. This proves the existence of three arcs $\Gamma_L$,
$\Gamma_C$, and $\Gamma_R$ contained in $\Gamma$ and connecting
$\zeta_{\pm}$, where $\xi_s \in \Gamma_s$ for $s\in \{L,C,R\}$.

Next, we note that $\Gamma_L \cup \Gamma_C$ is a closed contour consisting
of trajectories. It follows from \cite[Lemma 8.3]{Pommerenke}
that it has to surround a pole. Similarly
$\Gamma_C \cup \Gamma_R$ has to surround a pole.
This can only happen if $\xi_L < -1 < \xi_C < 1 < \xi_R$.

To complete the proof of the lemma, we need to establish
that $\Gamma$ consists only of $\Gamma_L$, $\Gamma_C$, and $\Gamma_R$
and nothing more.
We use that the function
\begin{equation} \label{defhz}
    h(z) = \Re \int_{\zeta_-}^z
    \frac{R(t)}{t^2-1} dt, \qquad
    z \in \mathbb C \setminus (\Gamma_C \cup \{-1,1\}),
\end{equation}
is single-valued and harmonic in $\mathbb C \setminus
(\Gamma_C \cup \{-1,1\})$.
The path of integration in (\ref{defhz}) is in $\mathbb C \setminus \Gamma_C$.
It is easy to see that
\[ \lim_{z \to \infty} h(z) = +\infty. \]
Since $h = 0$ on $\Gamma_L \cup \Gamma_R = \partial \Omega_{\infty}$,
it follows by the maximum principle for harmonic functions
that $h(z) > 0$ for $z \in \Omega_{\infty}$.
Similarly, since
\[ \lim_{z \to \infty} h(z) = - \infty, \]
and $h = 0$ on $\Gamma_L \cup \Gamma_C = \partial \Omega_{-1}$,
and on $\Gamma_R \cup \Gamma_C = \partial \Omega_1$,
we have that $h(z) < 0$ for $z \in \Omega_{\pm}$.
Since $\Gamma = \{h = 0\}$, we get that $\Gamma$
consists exactly of $\Gamma_L$, $\Gamma_C$, and $\Gamma_R$.
This completes the proof of Lemma \ref{lemma1}.
\end{varproof}

\begin{varproof} \textbf{of Lemma \ref{lemma2}.}
Recall that $ R(z):=\sqrt{(z-\zeta_+)(z-\zeta_-)}$ denotes
the single-valued branch in $\C \setminus \Gamma_C$ such that
$R(z) \sim z$ as $z \to \infty$. With this convention and taking
into account (\ref{zetapm}) it is straightforward to check that
\begin{equation}\label{R_1}
    R(-1)= \frac{2B}{A+B+2} <0, \qquad R(1)=\frac{-2A}{A+B+2} > 0.
\end{equation}
From the definition of $\Gamma$ it follows that $d \mu(z)$ is
real-valued on $\Gamma$ and does not change sign on each component
of $\Gamma \setminus \{\zeta_-,\zeta_+ \}$.

Using the residue theorem, we have that
\[
\mu(\Gamma_C \cup \Gamma_R)=\int_{\Gamma_C \cup \Gamma_R} d\mu
(t)= (A + B+2) \res_{z=1} \left( \frac{R(z)}{z^2-1} \right) = -A
\]
where we have used (\ref{R_1}).
Analogously,
\[
\mu(\Gamma_L \cup \Gamma_C)=\int_{\Gamma_L \cup \Gamma_C} d\mu
(t)= (A + B+2) \res_{z=-1} \left( \frac{R(z)}{z^2-1} \right) =-B.
\]
Finally,
\[
\mu(\Gamma_L \cup \Gamma_R)=\int_{\Gamma_L \cup \Gamma_R} d\mu
(t)= (A + B+2) \res_{z=\infty} \left( \frac{R(z)}{z^2-1} \right)
=A+B+2.
\]
Hence,
\[
\mu(\Gamma)=\frac{1}{2}\, \left( \mu(\Gamma_C \cup
\Gamma_R)+\mu(\Gamma_L \cup \Gamma_C) + \mu(\Gamma_L \cup
\Gamma_R)\right)=1,
\]
and
\begin{eqnarray*}
\mu(\Gamma_L) & = & 1-\mu(\Gamma_C \cup \Gamma_R)\ =\ 1+A, \\
\mu(\Gamma_R) & = & 1-\mu(\Gamma_L \cup \Gamma_C)\ = \ 1+B, \\
\mu(\Gamma_C) & = & 1 - \mu(\Gamma_L \cup \Gamma_R)\ =\ -1-A-B,
\end{eqnarray*}
which proves (\ref{measures_arcs1}). Since each part has positive total
$\mu$-mass and $\mu$ does not change sign on each of the parts,
we find that $\mu$ is a positive measure. This completes the proof.
\end{varproof}

\section{A Riemann-Hilbert problem for Jacobi polynomials}

\begin{figure}[htb]
\centering \includegraphics[scale=0.6]{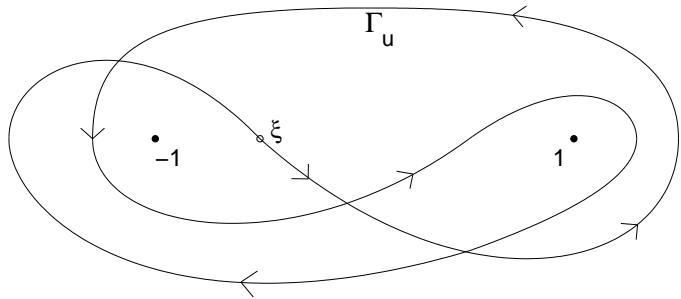} \caption{
The universal curve $\Gamma_u$.}\label{fig:pathGamma}
\end{figure}

Consider a closed path $\Gamma_u$ encircling the points $+1$ and $-1$
first in the positive direction and then in the negative
direction, as shown in Fig.~\ref{fig:pathGamma}. The point $\xi \in (-1,1)$ is the
begin and endpoint of $\Gamma_u$.

For $\alpha , \beta  \in \C$, denote
\begin{equation*}\label{weight}
w(z; \alpha , \beta ):=(1-z)^\alpha (1+z)^\beta=\exp [\alpha \log
(1-z) + \beta \log (1+z)].
\end{equation*}
This is a multi-valued function with branch points at $\infty$ and
$\pm 1$. However, if we start with a value of $w(z;\alpha, \beta)$ at a particular
point of $\Gamma_u$, and extend the definition of $w(z; \alpha, \beta)$
continuously along $\Gamma_u$, then we obtain a single-valued function
$w(z; \alpha, \beta)$ on $\Gamma_u$ if we view $\Gamma_u$ as a contour on
the Riemann surface for the function $w(z;\alpha, \beta)$.
For definiteness,
we assume that  the ``starting point'' is $\xi \in (-1,1)$, and
that the branch of $w$ is such that $w(\xi; \alpha , \beta )>0$.
We prefer to view $\Gamma_u$ as a subset of the complex plane.
Then $\Gamma_u$ has points of self-intersection, as shown
in  Fig.~\ref{fig:pathGamma}.
At points of self-intersection the value of $w(z;\alpha,\beta)$ is not
well-defined.

In \cite{Kuijlaars03} it was shown that for $k \in \{0, 1, \ldots, n\}$,
we have
\begin{equation}\label{orthMain}
    \int_{\Gamma_u} t^k \, P_n^{(\alpha , \beta )}(t) w(t; \alpha,
    \beta)\,
    dt = \frac{-  \pi^2 2^{n+\alpha +\beta +3} e^{\pi i (\alpha +\beta
    )}}{\Gamma(2n+\alpha +\beta +2) \Gamma(-n-\alpha )
    \Gamma(-n-\beta)}\, \delta_{kn}.
\end{equation}
This shows that the Jacobi polynomials $P_n^{(\alpha,\beta)}$
are orthogonal on the universal curve $\Gamma_u$.
The right-hand side of (\ref{orthMain}) vanishes for
$k=n$ if and only if either $-2n-\alpha-\beta-2$, or $n+\alpha$ or
$n+\beta$ is a non-negative integer. In some of these cases the
zero comes from integrating a single-valued and analytic function
along a curve in the region of analyticity; other values of
$\alpha$ and $\beta$ correspond to the special cases mentioned
before when there is a zero at $\pm 1$. It is shown in
\cite{Kuijlaars03} that the orthogonality conditions (\ref{orthMain})
characterize the Jacobi polynomial $P_n^{(\alpha, \beta)}$
provided the parameters satisfy
\begin{equation}\label{condition}
-n-\alpha-\beta \notin \N, \quad \text{and} \quad n+\alpha
    \notin \N, \quad \text{and} \quad n+\beta \notin \N.
\end{equation}
Then $P_n^{(\alpha , \beta )}$ is of degree exactly $n$, and we
will denote by $\widehat P_n^{(\alpha , \beta )}$ the
corresponding monic Jacobi polynomial.

Based on the orthogonality (\ref{orthMain})  a Riemann-Hilbert
problem is constructed in \cite{Kuijlaars03}, whose solution is
given in terms of  $\widehat P_n^{(\alpha,\beta)}$ with parameters
satisfying (\ref{condition}).

Let $\Gamma_u$ be a curve in $\C$ as described above with three
points of self-intersection as in Fig.\ \ref{fig:pathGamma}. We
let $\Gamma_u^o$ be the curve without the points of
self-intersection. Recall that the orientation of $\Gamma_u$ (see
also Fig.\ \ref{fig:pathGamma}) induces a $+$ and $-$ side in a
neighborhood of $\Gamma_u$, where the $+$ side is on the left
while traversing $\Gamma_u$ according to its orientation and the
$-$ side is on the right. Again, we say that a function $Y$ on $\C
\setminus \Gamma_u$ has a boundary value $Y_+(t)$ for $t \in
\Gamma_u^o$ if the limit of $Y(z)$ as $z \to t$ with $z$ on the
$+$ side of $\Gamma_u$ exists; similarly for $Y_-(t)$.

The Riemann-Hilbert problem for Jacobi polynomials is then as
follows. We look for a $2\times 2$ matrix valued function $Y
=Y^{(\alpha , \beta )}:\, \C \setminus \Gamma_u \to \C^{2\times
2}$ such that the following four conditions are satisfied:
\begin{enumerate}
\item[(a)] $Y$ is analytic on $\C \setminus \Gamma_u$.
\item[(b)]
$Y$ has continuous boundary values on $\Gamma_u^o$, denoted by $Y_{+}$
and $Y_{-}$, such that
\[ Y_{+}(t) = Y_{-}(t)
    \begin{pmatrix} 1 & w(t; \alpha, \beta) \\
        0 & 1 \end{pmatrix}
        \qquad \mbox{for } t\in \Gamma_u^o. \]
\item[(c)] As $z\to\infty$,
\[ Y(z) = \left(I + O\left(\frac{1}{z}\right)\right)
    \begin{pmatrix} z^{n} & 0 \\ 0 & z^{-n} \end{pmatrix}. \]
\item[(d)]
$Y(z)$ remains bounded as $z \to t \in \Gamma_u \setminus \Gamma_u^o$.
\end{enumerate}

This Riemann-Hilbert problem is similar to the Riemann-Hilbert
problem for orthogonal polynomials due to Fokas, Its, and Kitaev
\cite{Fokas92}, see also \cite{MR2000g:47048}.
Also the solution is similar.

\begin{proposition}[\cite{Kuijlaars03}] \label{solRH}
Assume that the parameters $\alpha, \beta$ satisfy
{\rm (\ref{condition})}. Then the above Riemann-Hilbert problem for
$Y$ has a unique solution, which is given by
\begin{equation} \label{formulaY}
    Y(z)= \begin{pmatrix}  \widehat P_{n}^{(\alpha,\beta)}(z) &
    \frac{1}{2\pi i} \int\limits_{\Gamma_u} \frac{\widehat P_{n}^{(\alpha,\beta)}(t)
        w(t; \alpha, \beta)}{t-z} \, dt \\[10pt]
    c_{n-1} P_{n-1}^{(\alpha,\beta)}(z) &
    \frac{c_{n-1}}{2\pi i} \int\limits_{\Gamma_u}
    \frac{P_{n-1}^{(\alpha,\beta)}(t) w(t; \alpha, \beta)}{t-z}\,  dt
    \end{pmatrix}, \qquad z \in \mathbb C \setminus \Gamma_u,
\end{equation}
for some non-zero constant $c_{n-1}$.
\end{proposition}

\medskip
The Riemann-Hilbert problem holds for any combination of
parameters $\alpha$ and $\beta$ such that (\ref{condition})
is satisfied.
Also the contour $\Gamma_u$ is rather arbitrary. It could be modified to
any curve that is homotopic to it in $\C \setminus \{-1,1\}$.

\section{Transformations of the Riemann-Hilbert problem}

In this section we consider parameters $A$ and $B$ satisfying the
inequalities (\ref{ABinequalities}). We also assume that $n \in
\mathbb N$ is such that $An$, $Bn$ and $(A+B)n$ are non-integers.
Throughout this section $A$, $B$, and $n$ remain fixed.

From Proposition \ref{solRH} we know that the Jacobi polynomial
$\widehat P_n^{(An,Bn)}$ is characterized as the $(1,1)$ entry of the
solution of the Riemann-Hilbert problem for $Y$ given in the previous
section with $\alpha = An$ and $\beta = Bn$. In this section we apply
the steepest descent method
of Deift and Zhou to this Riemann-Hilbert problem in order to
reduce it to a Riemann-Hilbert problem that is normalized at infinity
and whose jump matrices are close to the identity. In the next section we
derive the asymptotic results from this analysis. The Deift/Zhou steepest
descent method proceeds through a number of transformations of
the original Riemann-Hilbert problem.

\subsection{Choice of contour}
\begin{figure}[htb]
\centering \includegraphics[scale=0.65]{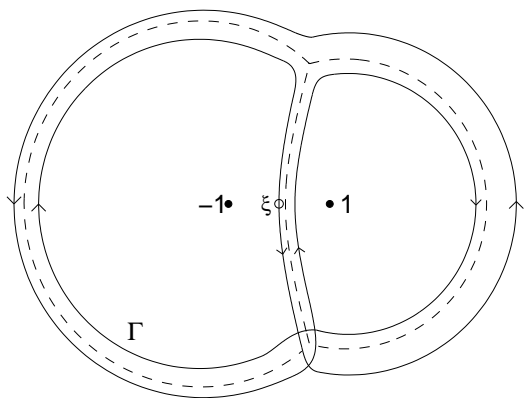}
\caption{Tautening $\Gamma_u$ on the set $\Gamma$.}
\label{fig_tautening}
\end{figure}

In the first step of the analysis we have to pick the right contour.
For  $A$ and $B$ satisfying (\ref{ABinequalities}) we have the
contour $\Gamma = \Gamma^{(A,B)}$ defined in (\ref{defGamma}),
which according to Lemma \ref{lemma1} consists of three analytic
arcs $\Gamma = \Gamma_L \cup \Gamma_C \cup \Gamma_R$.
We modify $\Gamma_u$ to a contour that follows $\Gamma$
in such a way that every part of $\Gamma$ is covered twice
as shown in Fig.\ \ref{fig_tautening}.

Passing from the  Riemann-Hilbert problem on $\Gamma_u$ to
the Riemann-Hilbert problem on $\Gamma$, we have that on each part
of $\Gamma$ two of the jumps are combined.
The new jump matrices take the form
\[ \begin{pmatrix} 1 & w(t_1; An, Bn) \\ 0 & 1 \end{pmatrix}
    \begin{pmatrix} 1 & -w(t_2; An, Bn) \\ 0 & 1 \end{pmatrix}
    = \begin{pmatrix} 1 & w(t_1; An,Bn) - w(t_2; An, Bn)
    \\ 0 & 1 \end{pmatrix},
    \]
where $t_1$ and $t_2$ are points on the Riemann surface, both
lying above $t$. The values of $w(t_j; An, Bn)$, $j=1, 2$,
differ from each other by a phase factor. To make this precise we
specify a single-valued branch for the weight
\[ w(z;An, Bn) = (1-z)^{An} (1+z)^{Bn} \]
on $\Gamma$. Since $\Gamma \setminus \{ \zeta_+\}$ is simply
connected, we can define a single valued branch on
$\Gamma \setminus \{ \zeta_+\}$, and we will do it in
such a way that $w(\xi; An, Bn) > 0$, where $\xi = \xi_C$ is
the intersection point of $\Gamma_C$ with the interval $(-1,1)$.

Then the jump on each of the contours $\Gamma_L$, $\Gamma_C$, and
$\Gamma_R$ can be calculated. The result is the following
Riemann-Hilbert problem on $\Gamma$ for a matrix valued function which we
continue to call $Y$.
The contour $\Gamma$ has the orientation shown in Fig.\ \ref{fig_Gamma1}.
\begin{enumerate}
\item[(a)] $Y$ is analytic on $\C \setminus \Gamma$.
\item[(b)] $Y$
has continuous boundary values on $\Gamma
\setminus \{\zeta_-, \zeta_+\}$, denoted by $Y_+$
and $Y_-$, such that
\begin{equation}\label{RHmatrix1}
Y_+(t) = Y_-(t)
    \begin{pmatrix} 1 & d_s w(t; An, Bn) \\
        0 & 1 \end{pmatrix}
        \qquad \mbox{for } t \in \Gamma_s \setminus \{ \zeta_{\pm}\},
            \quad s \in \{L,C,R\},
\end{equation}
with constants
\begin{equation}\label{constants}
d_L =  e^{2\pi Bn i} \left(e^{2\pi An i} -1\right), \qquad
   d_C =  1-e^{2\pi (A+B)ni}, \qquad
   d_R =  1-e^{2\pi B n i},
\end{equation}
and we follow the convention about the branch of $w(t; An, Bn)$
on $\Gamma$ mentioned above.
\item[(c)] As $z\to\infty$,
\[ Y(z) = \left(I + O\left(\frac{1}{z}\right)\right)
    \begin{pmatrix} z^{n} & 0 \\ 0 & z^{-n} \end{pmatrix}. \]
\item[(d)] $Y(z)$ remains bounded as $z \to \zeta_\pm $.
\end{enumerate}

Of course the solution to the above Riemann-Hilbert problem is similar
to the solution (\ref{formulaY}) to the earlier Riemann-Hilbert problem.
In particular we still have
\begin{equation} \label{Y11}
Y_{11}(z) = \widehat P_n^{(An,Bn)}(z)
\end{equation}

The constants $d_L$, $d_C$ and $d_R$ from (\ref{constants}) will play an important
role in what follows. These numbers are non-zero, exactly because of our
assumption that $An$, $Bn$, and $(A+B)n$ are non-integers. Observe that
\begin{equation}\label{constants2}
d_L + d_C =  d_R,
\end{equation}
which is a relation that will be used a number of times.

\subsection{Auxiliary functions}
In order to make the first transformation of the Riemann-Hilbert problem
we need some auxiliary functions.

We already know from Lemma \ref{lemma2} that $\mu$
defined in (\ref{defmu}) is a positive measure on $\Gamma$ such
that
\begin{equation}\label{measures_arcs}
\mu(\Gamma_L)=1+A>0, \qquad \mu(\Gamma_C)=-1-A-B>0, \qquad
\mu(\Gamma_R)=1+B>0.
\end{equation}

Let $g$ be the complex logarithmic potential of the measure $\mu$,
\[ g(z) = \int \log(z-t) d\mu(t). \]
This is a multivalued function; however its derivative is single
valued:
\begin{equation*}\label{g_prime}
g'(z)=\int \frac{d\mu(t)}{z-t}=\begin{cases} \dfrac{\strut
A+B+2}{2}\, \dfrac{R(z)}{z^2-1}-
\dfrac{A/2}{z-1}-\dfrac{B/2}{z+1}, & \text{ for } z
    \in \Omega_\infty,\\[10pt]
    -\dfrac{\strut A+B+2}{2}\, \dfrac{R(z)}{z^2-1}-
\dfrac{A/2}{z-1}-\dfrac{B/2}{z+1}, & \text{ for } z
    \in \Omega_{-1} \cup \Omega_{1}.
\end{cases}
\end{equation*}

We define
\begin{equation*}\label{def_G}
    G(z)=\exp \left(\int_{\zeta_-}^z g'(t)\, dt\right), \qquad  z \in \mathbb C \setminus \Gamma,
\end{equation*}
where the path of integration lies entirely in $\mathbb C
\setminus \Gamma$ except for the initial point $\zeta_-$. From the
fact that $\mu$ is a positive unit measure on $\Gamma$ it follows
that $G$ is single-valued in each component of $\mathbb C
\setminus \Gamma$. Furthermore, $G$ is analytic, $G(\zeta_-) = 1$,
and the following limit exists
 \begin{equation}\label{kappa}
   \kappa := \lim_{z \to \infty} \frac{G(z)}{z} = \zeta_-
    \exp\left( \int_{\zeta_-}^{\infty} (g'(t) - 1/t)\, dt  \right).
\end{equation}
We calculate the jumps of $G$.
We have
\begin{equation} \label{G+timesG-}
    G_+(z) G_-(z) = \frac{w(\zeta_-;A,B)}{w(z;A,B)}, \qquad \textrm{ for } z \in \Gamma,
\end{equation}
and
\begin{equation} \label{G+divideG-}
    \frac{G_+(z)}{G_-(z)} = \exp(-2 \phi_+(z)), \qquad  \textrm{ for } z \in \Gamma,
\end{equation}
where $\phi$ is defined by (\ref{defphi2}).
It will be useful to introduce also the related function
\begin{equation} \label{defphitilde}
    \tilde{\phi}(z) = \frac{A+B+2}{2} \int_{\zeta_+}^z \frac{R(t)}{t^2-1} dt
    = \overline{\phi(\bar{z})},
        \qquad \text{ for } z \in \mathbb C \setminus (\Gamma_C \cup \gamma_{-1}^{-}
            \cup \gamma_1^- \cup \gamma_{\infty}^-).
\end{equation}
To relate $\tilde{\phi}$ with $\phi$ it is necessary to compute
$\frac{A+B+2}{2} \int_{\zeta_-}^{\zeta_+} \frac{R(t)}{t^2-1} dt$. This
integral depends on the path from $\zeta_-$ to $\zeta_+$. We
can distinguish four paths, namely $\Gamma_R$, $-\Gamma_{C,+}$,
$-\Gamma_{C,-}$ and $-\Gamma_L$. (Recall that $\Gamma_C$ and
$\Gamma_L$ are oriented from $\zeta_+$ to $\zeta_-$. So we put a
minus sign to indicate that the path is from $\zeta_-$ to
$\zeta_+$.) We obtain
\[ \frac{A+B+2}{2} \int_{\zeta_-}^{\zeta_+} \frac{R(t)}{t^2-1} dt
    = \left\{ \begin{array}{lcll}
    \pi i \mu(\Gamma_R) & = & \pi i (1+B) & \text{ integral over $\Gamma_R$} \\
    -\pi i \mu(\Gamma_C) & = & \pi i (1+A+B) & \text{ integral over $-\Gamma_{C,+}$} \\
    \pi i \mu(\Gamma_C) & = & -\pi i (1+A+B) & \text{ integral over $-\Gamma_{C,-}$} \\
    -\pi i \mu(\Gamma_L) & = & -\pi i (1+A) & \text{ integral over $-\Gamma_L$}
    \end{array} \right. \]
where we have used (\ref{measures_arcs}). It follows that
\begin{eqnarray} \label{phirelation1}
   \left. \begin{array}{rcl}
     \phi_+(z) & = & \tilde{\phi}(z) + \pi i(1+B) \qquad \\
     \phi_-(z) & = & \tilde{\phi}(z) - \pi i(1+A)
     \end{array} \right\} && \text{ for $z$ on } \gamma_{\infty}^+, \\
     \label{phirelation2}
   \left. \begin{array}{rcl}
     \phi_+(z) & = & \tilde{\phi}(z) - \pi i(1+A+B) \\
     \phi_-(z) & = & \tilde{\phi}(z) - \pi i(1+A)
     \end{array} \right\} && \text{ for $z$ on } \gamma_{-1}^+, \\
     \label{phirelation3}
   \left. \begin{array}{rcl}
     \phi_+(z) & = & \tilde{\phi}(z) + \pi i(1+B) \\
     \phi_-(z) & = & \tilde{\phi}(z) + \pi i(1+A+B)
     \end{array} \right\} && \text{ for $z$ on } \gamma_{1}^+.
\end{eqnarray}
Observe also that by construction both $\phi$ and $\tilde{\phi}$
have negative real parts in the bounded components $\Omega_{-1}$
and $\Omega_{1}$ of
$\C \setminus \Gamma$ (where defined) and positive real part in
$\Omega_\infty$ (with the appropriate cuts).

\subsection{First transformations $Y \mapsto U$}

Now we introduce a new matrix valued function $U$ by
\begin{equation} \label{defUz}
    U(z) = \kappa^{n\sigma_3}\, w(\zeta_-; An ,Bn )^{-\sigma_3/2}\,
    Y(z)
    \, G(z)^{-n\sigma_3}\, w(\zeta_-; An ,Bn )^{\sigma_3/2},
\end{equation}
where $\sigma_3 =\begin{pmatrix} 1 & 0 \\ 0& -1 \end{pmatrix}$ is
the Pauli matrix, and for any non-zero $x$, $
x^{\sigma_3}=\begin{pmatrix} x & 0 \\
0& 1/x
\end{pmatrix}
$. Here $G$ is the function introduced in (\ref{def_G}), and
$\kappa$ is the limit defined in (\ref{kappa}). Then $U$ satisfies
the Riemann-Hilbert problem
\begin{enumerate}
\item[(a)] $U$ is analytic on $\C \setminus \Gamma$.
\item[(b)]
$U$ has continuous boundary values on $\Gamma \setminus \{\zeta_\pm\}$
such that
\begin{equation} \label{jumpU}
    U_{+}(z) = U_{-}(z)
    \begin{pmatrix} \exp(2n \phi_+(z)) & d_s \\ 0 & \exp(-2n \phi_+(z))
    \end{pmatrix}
        \qquad \mbox{for } z\in \Gamma_s^o,
            \quad s \in \{L,C,R\}.
\end{equation}
\item[(c)] $ U(z) = I + O\left(\frac{1}{z}\right)$ as $z \to
\infty$. \item[(d)] $U(z)$ remains bounded as $z \to \zeta_{\pm}$.
\end{enumerate}

To obtain the jumps in (\ref{jumpU}) we used the relations (\ref{G+timesG-}) and
(\ref{G+divideG-}). For the asymptotic behavior in (c) we used the limit
(\ref{kappa}).

We use the following factorizations of the jump matrices in (\ref{jumpU})
\begin{equation} \label{factorjump1}
    \begin{pmatrix} e^{2n\phi_+} & d_C \\ 0 & e^{-2n\phi_+} \end{pmatrix}
   = \begin{pmatrix} 1 & 0 \\ \frac{1}{d_C} e^{2n\phi_-} & 1 \end{pmatrix}
     \begin{pmatrix} 0 & d_C \\ - \frac{1}{d_C} & 0 \end{pmatrix}
     \begin{pmatrix} 1 & 0 \\ \frac{1}{d_C} e^{2n\phi_+} & 1 \end{pmatrix}
\end{equation}
and
\begin{equation} \label{factorjump2}
    \begin{pmatrix} e^{2n\phi} & d_s \\ 0 & e^{-2n\phi} \end{pmatrix}
   = \begin{pmatrix} 1 & 0 \\ \frac{1}{d_s} e^{-2n\phi} & 1 \end{pmatrix}
     \begin{pmatrix} e^{2n\phi} & d_s \\ -\frac{1}{d_s} & 0 \end{pmatrix},
     \quad \text{ for } s=L,R,
 \end{equation}
in order to define the next transformation.

\subsection{Second transformation $U \mapsto T$}

The trajectories $\Gamma$ and the orthogonal trajectories
$\Gamma^{\perp}$ divide the complex plane into six domains, which
we number from left to right as domains I, II, III, IV,
V and VI, see Fig.\ \ref{fig_regions1}. We define
$\tilde{T}$ in each of these six
domains separately. We put
\begin{align}
\tilde{T} &= U \label{tildeT1}
    \begin{pmatrix} 1 & 0 \\ \frac{1}{d_L} e^{-2n\phi} & 1 \end{pmatrix}
    \quad \text{in domain I}, \\
\tilde{T} &= U \label{tildeT2}
    \begin{pmatrix} 1 & 0 \\ \frac{1}{d_R} e^{-2n\phi} & 1 \end{pmatrix}
    \quad \text{in domain VI}, \\
\tilde{T} &=U \label{tildeT3}
    \begin{pmatrix} 1 & 0  \\ -\frac{1}{d_L} e^{2n\phi} & 1 \end{pmatrix}
    \begin{pmatrix} 0 & -d_L \\ \frac{1}{d_L} & 0 \end{pmatrix}= U
    \begin{pmatrix} 0 & -d_L \\ \frac{1}{d_L} & e^{2n\phi} \end{pmatrix}
    \quad \mbox{in domain II,} \\
\tilde{T} &=U \label{tildeT4}
    \begin{pmatrix} 1 & 0  \\ -\frac{1}{d_R} e^{2n\phi} & 1 \end{pmatrix}
    \begin{pmatrix} 0 & -d_R \\ \frac{1}{d_R} & 0 \end{pmatrix}= U
    \begin{pmatrix} 0 & -d_R \\ \frac{1}{d_R} & e^{2n\phi} \end{pmatrix}
    \quad \text{in domain V,} \\
\tilde{T} &= U \label{tildeT5}
    \begin{pmatrix} 1 & 0  \\ \frac{1}{d_C} e^{2n\phi} & 1 \end{pmatrix}
    \begin{pmatrix} 0 & -d_L \\ \frac{1}{d_L} & 0 \end{pmatrix}
    \quad \mbox{in domain III,} \\
\tilde{T} & = U \label{tildeT6}
    \begin{pmatrix} 1 & 0  \\ -\frac{1}{d_C} e^{2n\phi} & 1 \end{pmatrix}
    \begin{pmatrix} 0 & -d_R \\ \frac{1}{d_R} & 0 \end{pmatrix}
    \quad \mbox{in domain IV}.
\end{align}
Since $\phi(z)$ behaves like $\frac{A+B+2}{2} \log z$ as $z \to
\infty$, we have $|e^{-2n \phi(z)}| \sim |z|^{-(A  +B  +2 )n}$, so
that
\[  \begin{pmatrix} 1 & 0 \\ \frac{1}{d_s}\, e^{-2n\phi(z)} & 1 \end{pmatrix}
    = I + O(1/z) \quad \mbox{ as } z \to \infty. \]
Thus $\tilde{T}(z) = I  + O(1/z)$ as $z \to \infty$.

By definition, $\tilde{T}$ is analytic in $\mathbb C \setminus
(\Gamma \cup \Gamma^{\perp})$. However, we have arranged our
transformation in a way that the jumps on $\Gamma_L$ and
$\Gamma_R$ disappear (due to the factorization (\ref{factorjump2})
and the definition of $\tilde{T}$) so $\tilde{T}$ is analytic in
$\mathbb C \setminus (\Gamma_C \cup \Gamma^{\perp})$.

We compute the jumps on $\Gamma_C \cup \Gamma^{\perp}$ with the convention that
these curves are oriented as shown in Fig.\ \ref{fig_regions1}. Straightforward
computations then show that
\begin{align*}
\tilde{T}_+ &= \tilde{T}_-
    \begin{pmatrix} 0 & \frac{d_Ld_R}{d_C} \\ - \frac{d_C}{d_Ld_R} & 0 \end{pmatrix}
        \quad \mbox{ on } \Gamma_C,\\
\tilde{T}_+ &= \tilde{T}_-
    \begin{pmatrix} 1 & 0 \\ \frac{1}{d_R} e^{-2n\phi_+} - \frac{1}{d_L} e^{-2n\phi_-} & 1 \end{pmatrix}
        \quad \mbox{ on } \gamma_{\infty}^+ \cup \gamma_{\infty}^-, \\
        \tilde{T}_+ &= \tilde{T}_-
    \begin{pmatrix} 1 & - d_L e^{2n\phi_-} - \frac{d_L^2}{d_C} e^{2n\phi_+} \\ 0 & 1 \end{pmatrix}
        \quad \mbox{ on } \gamma_{-1}^+ \cup \gamma_{-1}^-,\\
        \tilde{T}_+ &= \tilde{T}_-
    \begin{pmatrix} 1 &  -\frac{d_R^2}{d_C} e^{2n\phi_-} + d_R e^{2n\phi_+} \\ 0 & 1 \end{pmatrix}
        \quad \mbox{ on } \gamma_{1}^+ \cup \gamma_{1}^-.
\end{align*}
Since $\phi$ is analytic across the curves $\gamma_{j}^-$, the jumps on these
curves simplify to (we also use (\ref{constants2}))
\begin{align*}
\tilde{T}_+ &= \tilde{T}_-
    \begin{pmatrix} 1 & 0 \\ - \frac{d_C}{d_Ld_R} e^{-2n\phi} & 1 \end{pmatrix}
        \quad \mbox{ on } \gamma_{\infty}^- , \\
        \tilde{T}_+ &= \tilde{T}_-
    \begin{pmatrix} 1 & -\frac{d_L d_R}{d_C} e^{2n\phi} \\ 0 & 1 \end{pmatrix}
        \quad \mbox{ on } \gamma_{-1}^- \cup \gamma_{1}^-.
\end{align*}
 If we now express the jumps on the contours $\gamma_{j}^+$ in terms of
$\tilde{\phi}$, see (\ref{phirelation1})--(\ref{phirelation3}),
they look as those on the lower half plane, but
with $\phi$ replaced by $\tilde \phi$:
\begin{align*}
\tilde{T}_+ &= \tilde{T}_-
    \begin{pmatrix} 1 & 0 \\ - \frac{d_C}{d_Rd_L} e^{-2n\tilde \phi} & 1 \end{pmatrix}
        \quad \mbox{ on } \gamma_{\infty}^+ , \\
        \tilde{T}_+ &= \tilde{T}_-
    \begin{pmatrix} 1 & -\frac{d_L d_R}{d_C} e^{2n\tilde \phi} \\ 0 & 1 \end{pmatrix}
        \quad \mbox{ on } \gamma_{-1}^+ \cup \gamma_{1}^+.
\end{align*}

Now with $\tau$ such that
\begin{equation} \label{def tau}
    \tau^2 = \frac{d_L d_R}{d_C}
\end{equation}
we define $T$ by
\begin{equation} \label{defT}
    T = \begin{pmatrix} \tau^{-1} & 0 \\ 0 & \tau \end{pmatrix}
    \tilde{T} \begin{pmatrix} \tau & 0 \\ 0 & \tau^{-1} \end{pmatrix}.
\end{equation}
The effect on the jump matrices is that the $(1,2)$ entries are
multiplied by $\tau^{-2}$ and the $(2,1)$ entries are multiplied
by $\tau^2$. So $T$ satisfies the following Riemann-Hilbert
problem:
\begin{enumerate}
\item[(a)] $T$ is analytic on $\C \setminus ( \Gamma_C \cup \Gamma^\perp)$.
\item[(b)]
$T$ has continuous boundary values on $(\Gamma_C \cup \Gamma^\perp)
\setminus \{\zeta_\pm\}$ such that
\begin{align}
T_+ &= T_- \label{jumpT1}
    \begin{pmatrix} 0 & 1 \\ -1 & 0 \end{pmatrix}
        \quad \text{ on } \Gamma_C, \\
 T_+ & = T_- \label{jumpT2}
    \begin{pmatrix} 1 & 0 \\ - e^{-2n\phi} & 1 \end{pmatrix}
        \quad \text{ on } \gamma_{\infty}^- ,\\
T_+ & = T_- \label{jumpT3}
    \begin{pmatrix} 1 & - e^{2n\phi} \\ 0 & 1 \end{pmatrix}
        \quad \text{ on }  \gamma_{-1}^- \cup \gamma_{1}^- ,\\
 T_+ & = T_- \label{jumpT4}
    \begin{pmatrix} 1 & 0 \\ - e^{-2n\tilde{\phi}} & 1 \end{pmatrix}
        \quad \text{ on } \gamma_{\infty}^+ ,\\
 T_+ & = T_- \label{jumpT5}
    \begin{pmatrix} 1 & - e^{2n\tilde{\phi}} \\ 0 & 1 \end{pmatrix}
        \quad \text{ on } \gamma_{-1}^+ \cup \gamma_{1}^+.
        \end{align}
\item[(c)]
$T(z) = I + O\left(\frac{1}{z}\right)$ as $z \to \infty$.
\item[(d)]
$T(z)$ remains bounded as $z \to \zeta_{\pm}$.
\end{enumerate}

The problem for $T$ is by now relatively standard. However, compared with earlier
works, the triangularity of the jump matrices on the curves $\gamma_j^{\pm}$ is
reversed. The inverse transposed matrix $T^{-t}$ satisfies the jumps
\begin{align*}
T_+^{-t} &= T_-^{-t}
    \begin{pmatrix} 0 & 1 \\ -1 & 0 \end{pmatrix}
        \quad \text{ on } \Gamma_C, \\
 T_+^{-t} & = T_-^{-t}
    \begin{pmatrix} 1 &  e^{-2n\phi}\\ 0 & 1 \end{pmatrix}
        \quad \text{ on } \gamma_{\infty}^- ,\\
T_+^{-t} & = T_-^{-t}
    \begin{pmatrix} 1 & 0\\ e^{2n\phi}  & 1 \end{pmatrix}
        \quad \text{ on } \gamma_{-1}^- \cup \gamma_{1}^- ,\\
 T_+^{-t} & = T_-^{-t}
    \begin{pmatrix} 1 &  e^{-2n\tilde{\phi}}\\ 0 & 1 \end{pmatrix}
        \quad \text{ on } \gamma_{\infty}^+ ,\\
 T_+ ^{-t}& = T_-^{-t}
    \begin{pmatrix} 1 & 0 \\ e^{2n\tilde{\phi}}  & 1 \end{pmatrix}
        \quad \text{ on } \gamma_{-1}^+ \cup \gamma_{1}^+,
        \end{align*}
which are exactly of the form considered for example in
\cite{MR2001g:42050,Kuijlaars/Mclaughlin:01a}.

\subsection{Outside parametrix}

The jump matrices in (\ref{jumpT2})--(\ref{jumpT5}) are close to the
identity matrix if $n$ is large. Therefore we expect that the main term in the
asymptotic behavior
of $T$ is given by the solution $N$ to the following model Riemann-Hilbert problem:
\begin{enumerate}
\item[(a)] $N$ is analytic in $\mathbb C \setminus \Gamma_C$,
\item[(b)] $N_+ = N_- \begin{pmatrix} 0 & 1 \\ -1 & 0 \end{pmatrix}$
         on $\Gamma_C \setminus \{\zeta_\pm\}$,
\item[(c)] $N(z) = I + O(1/z)$ as $z \to \infty$.
\end{enumerate}
In analogy with the condition (d) in the Riemann-Hilbert problem for $T$
we would like to ask that $N(z)$ remains bounded as $z \to \zeta_{\pm}$.
However, this would lead to a Riemann-Hilbert problem with no solution.
Instead we allow for moderate singularities of $N$ at $\zeta_{\pm}$:
\begin{enumerate}
\item[(d)] $N(z)=O(|z-\zeta_{\pm}|^{-1/4})$ as $z \to \zeta_\pm$.
\end{enumerate}

The solution to this problem is given by
\begin{equation} \label{defN}
    N(z) = \begin{pmatrix} \frac{a(z)+a(z)^{-1}}{2} & \frac{a(z) - a(z)^{-1}}{2i} \\
    -\frac{a(z)-a(z)^{-1}}{2i} & \frac{a(z) + a(z)^{-1}}{2} \end{pmatrix}
\end{equation}
with
\[ a(z) = \frac{(z-\zeta_-)^{1/4}}{(z-\zeta_+)^{1/4}}, \qquad
    z \in \mathbb C \setminus \Gamma_C, \]
being analytic and single-valued in $\C \setminus \Gamma_C$,
such that $a(z) \to 1$ as $z \to \infty$, see
\cite{MR2000g:47048,MR2001f:42037,MR2001g:42050,Kuijlaars/Mclaughlin:01a}.
In
\cite{Kuijlaars/Mclaughlin:01a} also an alternative expression for $N$
has been established in terms of $R(z):=\sqrt{(z-\zeta_+)(z-\zeta_-)}$:
\begin{equation}\label{alternativeForN}
N(z) = \begin{pmatrix} \left(\frac{1+R'(z)}{2}\right)^{1/2} & -\left(\frac{1-R'(z)}{2}\right)^{1/2} \\
    \left(\frac{1-R'(z)}{2}\right)^{1/2} &
    \left(\frac{1+R'(z)}{2}\right)^{1/2}\end{pmatrix}.
\end{equation}

\subsection{Local parametrices}

Near the branch points $\zeta_{\pm}$ we construct local
parametrices in the same way as done by Deift et al
\cite{MR2001f:42037,MR2001g:42050,MR2000g:47048}, see also
\cite{Kuijlaars/Mclaughlin:01b,Kuijlaars/Mclaughlin:01a}. In a
neighborhood $U_{\delta} =\{ z \in \mathbb C :\, |z-\zeta_-| <
\delta\}$ of $\zeta_-$ we construct a $2\times 2$ matrix valued
$P$ that is analytic in $U_{\delta} \setminus (\Gamma_C \cup
\gamma_{-1}^- \cup \gamma_1^- \cup \gamma_{\infty}^-)$, satisfies
the same jump conditions as $T$ does on $U_{\delta} \cap (\Gamma_C
\cup \gamma_{-1}^- \cup \gamma_1^- \cup \gamma_{\infty}^-)$ and
that matches with $N$ on the boundary $C_{\delta}$ of $U_{\delta}$
up to order $1/n$.

\begin{figure}[htb]
\centering \includegraphics[scale=0.65]{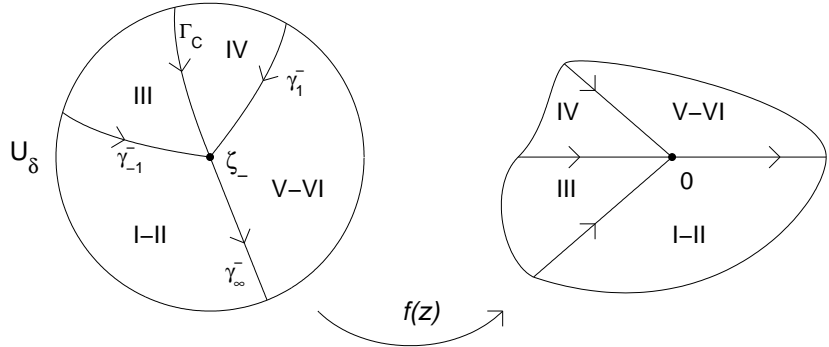}
\caption{Conformal mapping $f$. }\label{fig_Jacobiconfmap}
\end{figure}

We need the function
\begin{equation}\label{deff}
f(z) = \left[ \frac{3}{2} \phi(z) \right]^{2/3}
\end{equation}
where the $2/3$rd root is chosen which is real and positive on
$\gamma_{\infty}^-$. This is a conformal map from $U_{\delta}$ onto
a neighborhood of $0$ provided $\delta > 0$ is small enough.
We note that $\gamma_{\infty}^-$ is mapped
to the positive real axis, $\Gamma_C$ to (a part of) the negative
real axis. Recall that $\phi$ is real and negative on $\gamma_1^-$
and $\gamma_{-1}^-$ and we see that $\gamma_1^-$ is mapped to
$\arg w = 2\pi/3$ and $\gamma_{-1}^-$ to $\arg w = - 2\pi/3$
(Fig.\ \ref{fig_Jacobiconfmap}).
Then the Riemann-Hilbert problem for $P$ is solved by (cf.\
\cite{Kuijlaars/Mclaughlin:01a})
\begin{equation} \label{defP}
    P(z) = \left[ E(z) \Psi(n^{2/3} f(z)) e^{n\phi(z)\sigma_3} \right]^{-t},
\end{equation}
where
\begin{equation} \label{defE}
    E(z) = \sqrt{\pi} e^{\frac{\pi i}{6}} \left(\begin{array}{cc}
        1 & -1 \\ -i & -i \end{array} \right)
        \left(\frac{n^{1/6} f(z)^{1/4}}{a(z)} \right)^{\sigma_3},
\end{equation}
and $\Psi$ is  built out of the Airy function $\Ai$ and its
derivative $\Ai'$ as follows
\begin{equation} \label{defPsi1}
    \Psi(s) = \left\{ \begin{array}{ll}
    \left(\begin{array}{cc}
    \Ai(s) & \Ai(\omega^2 s) \\
    \Ai'(s) & \omega^2 \Ai'(\omega^2 s) \end{array} \right)
    e^{- \frac{\pi i}{6} \sigma_3} &
    \mbox{for } 0 < \arg s < 2 \pi/3, \\[10pt]
    \left(\begin{array}{cc}
    \Ai(s) & \Ai(\omega^2 s) \\
    \Ai'(s) & \omega^2 \Ai'(\omega^2 s) \end{array} \right)
    e^{- \frac{\pi i}{6} \sigma_3}
    \left(\begin{array}{cc} 1 & 0 \\ -1 & 1 \end{array} \right) &
    \mbox{for } 2\pi/3 < \arg s < \pi, \\[10pt]
    \left(\begin{array}{cc}
    \Ai(s) & -\omega^2 \Ai(\omega s) \\
    \Ai'(s) & -\Ai'(\omega s) \end{array} \right)
    e^{-\frac{\pi i}{6} \sigma_3}
    \left(\begin{array}{cc} 1 & 0 \\ 1 & 1 \end{array} \right) &
    \mbox{for } -\pi < \arg s < -2\pi/3, \\[10pt]
    \left(\begin{array}{cc}
    \Ai(s) & -\omega^2 \Ai(\omega s) \\
    \Ai'(s) & -\Ai'(\omega s) \end{array} \right)
    e^{- \frac{\pi i}{6} \sigma_3} &
    \mbox{for } -2\pi/3 < \arg s < 0,
    \end{array} \right.
\end{equation}
with $\omega = e^{2\pi i/3}$.

Note that we take the inverse transpose in (\ref{defP}), which is
absent in the construction in \cite{Kuijlaars/Mclaughlin:01a}.  This is of course
due to the fact that the Riemann-Hilbert problem for $T^{-t}$ is
comparable to the Riemann-Hilbert problem found in
\cite{Kuijlaars/Mclaughlin:01a}, see the remark at the
end of subsection 5.4.
\medskip

A similar construction yields a parametrix $\tilde{P}$ in a neighborhood
$\tilde{U}_{\delta} = \{ z : |z-\zeta_+| < \delta\}$.

\subsection{Final transformation $T \mapsto S$}

\begin{figure}[htb]
\centering \includegraphics[scale=0.65]{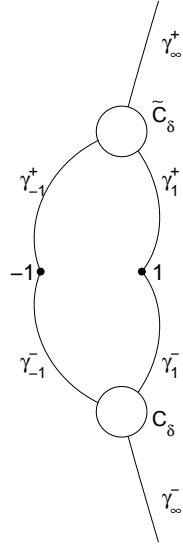}
\caption{Contour $\Gamma_S$ for the Riemann-Hilbert problem of
$S$.}\label{fig_Jac_fig36_finaljumps}
\end{figure}

The final transformation $T \mapsto S$ is
\begin{eqnarray}
\label{defS1}
S = TN^{-1} && \text{ outside the disks $U_{\delta}$ and $\tilde{U}_{\delta}$}, \\
\label{defS2}
S = TP^{-1} && \text{ inside the disk $U_{\delta}$}, \\
\label{defS3}
S = T \tilde{P}^{-1} && \text{ inside the disk $\tilde{U}_{\delta}$}.
\end{eqnarray}
Then by construction, $S$ has jumps on the circles $C_{\delta}=\partial
U_{\delta}$ and $\tilde{C}_{\delta}=\partial \tilde{U}_{\delta}$
as well as on $\Gamma_C \cup \Gamma^\perp$. Since the jumps of $T$ and $N$ on
$\Gamma_C$ agree, we have that $S$ is analytic
across the part of $\Gamma_C$ outside the disks $U_{\delta}$ and $\tilde{U}_{\delta}$.
Similarly, the jumps of $T$ and $P$ agree inside the disk $U_{\delta}$,
and the jumps of $T$ and $\tilde{P}$ agree inside the disk $\tilde{U}_{\delta}$,
so that $S$ is analytic in $U_{\delta}$ and $\tilde{U}_{\delta}$ with possible
isolated singularities at $\zeta_{\pm}$. However it follows from the behavior
of $T$ and $N$ near $\zeta_{\pm}$ that the singularities are
removable. Thus $S$ solves the following Riemann-Hilbert problem.

\begin{enumerate}
\item[(a)] $S$ is analytic on $\mathbb C \setminus \Gamma_S$, where $\Gamma_S$
consists of the circles $C_{\delta}$ and $\tilde{C}_{\delta}$, and of the
parts of $\gamma_{-1}$, $\gamma_1$ and $\gamma_{\infty}$ outside the disks,
see Fig.\ \ref{fig_Jac_fig36_finaljumps}.
\item[(b)]
$S$ has continuous boundary values on $\Gamma_S$ such that
\begin{align}
S_+ & = S_- \label{jumpS1}
    N \begin{pmatrix} 1 & 0 \\ - e^{-2n\phi} & 1 \end{pmatrix} N^{-1}
        \quad \text{ on } \gamma_{\infty}^- \setminus \tilde{U}_{\delta},\\
S_+ & = S_- \label{jumpS2}
    N \begin{pmatrix} 1 & - e^{2n\phi} \\ 0 & 1 \end{pmatrix} N^{-1}
        \quad \text{ on }  \left(\gamma_{-1}^- \cup \gamma_{1}^- \right)
            \setminus \tilde{U}_{\delta},\\
S_+ & = S_- \label{jumpS3}
    N \begin{pmatrix} 1 & 0 \\ - e^{-2n\tilde{\phi}} & 1 \end{pmatrix} N^{-1}
        \quad \text{ on } \gamma_{\infty}^+\setminus U_{\delta},\\
S_+ & = S_- \label{jumpS4}
    N \begin{pmatrix} 1 & - e^{2n\tilde{\phi}} \\ 0 & 1 \end{pmatrix} N^{-1}
        \quad \text{ on } \left( \gamma_{-1}^+ \cup \gamma_{1}^+\right)
            \setminus U_{\delta}, \\
S_+ & = S_- \label{jumpS5}
    P N^{-1} \qquad \qquad \text{ on } C_{\delta}, \\
S_+ & = S_- \label{jumpS6}
    \tilde{P} N^{-1} \qquad \qquad \text{ on } \tilde{C}_{\delta}.
\end{align}
\item[(c)] $S(z) = I + O\left(\frac{1}{z}\right) $ as $z \to \infty$.
\end{enumerate}

\section{Asymptotics: Proofs of the theorems}

\subsection{Asymptotics of $S$}

The analysis in the last section is done for fixed values of $A$, $B$,
and $n$. All the transformations are exact for finite $n$. It is now
our aim to let $n \to \infty$ and control the jump matrices in the
Riemann-Hilbert problem for $S$. We want to do it in a way which
is valid locally uniformly for parameters $A$ and $B$
satisfying (\ref{ABinequalities}).

Then first of all we should study the dependence of the contour
$\Gamma_S$ on the parameters $A$ and $B$. Note that $\Gamma_S$ does not
depend on $n$, but it does depend on $A$ and $B$. In fact, we have
that $\Gamma^\perp$ is completely determined by $A$ and $B$, while
the radius $\delta$ of the circles around $\zeta_{\pm}$
is only restricted by the requirement that the mapping $f$
from (\ref{deff}) is a conformal mapping on $U_{\delta}$.
With that in mind, it is clear that we may assume that the curve
$\Gamma_S$ depends on $A$ and $B$ in a continuous way.

Now we can see what happens with the jump matrices in
(\ref{jumpS1})--(\ref{jumpS6}) as $n \to \infty$.
On $\gamma_{\infty}^- \setminus \tilde{U}_{\delta}$ we have
that $\Re \phi$ is strictly positive.
Hence the jump matrix in (\ref{jumpS1}) is $I + O(e^{-cn})$ as $n \to \infty$,
uniformly on $\gamma_{\infty}^- \setminus \tilde{U}_{\delta}$.
This estimate is also valid uniformly for $A$ and $B$ in compact subsets
of the set
\begin{equation} \label{ABset}
    \{ (A,B) \mid -1 < A < 0,\ -1 < B < 0,\ -2 < A+B< -1 \}.
\end{equation}
Similarly, the jump matrices in
(\ref{jumpS2})--(\ref{jumpS4}) are $I + O(e^{-cn})$ as $n \to \infty$,
uniformly on the respective contours and uniformly for $A$ and $B$
in compact subsets of (\ref{ABset}).

For (\ref{jumpS5}) and (\ref{jumpS6}) we make use of the matching conditions
\begin{equation} \label{matching1}
    P(z) = \left(I + O\left(\frac{1}{n}\right)\right) N(z)
    \qquad \mbox{uniformly for } z \in C_{\delta}.
\end{equation}
and
\begin{equation} \label{matching2}
\tilde{P}(z) = \left(I + O\left(\frac{1}{n}\right)\right) N(z)
    \qquad \mbox{uniformly for } z \in \tilde{C}_{\delta}.
\end{equation}
So that the jump matrices in (\ref{jumpS5}) and (\ref{jumpS6}) are $I + O(1/n)$
as $n \to \infty$, uniformly on the two circles.
A closer analysis also reveals that the $O$-terms in (\ref{matching1})
and (\ref{matching2}) are valid uniformly for $A$ and $B$ in compact
subsets of (\ref{ABset}).

So the conclusion is that all jumps in (\ref{jumpS1})--(\ref{jumpS6}) are
$I + O(1/n)$ uniformly for $z$ on $\Gamma_S$, and uniformly for $A$ and $B$
in compact subsets of (\ref{ABset}).
Then arguments such as in \cite{MR2000g:47048,MR2001f:42037,MR2001g:42050}
lead to the following conclusion.
\begin{proposition}
We have that
\begin{equation} \label{asympS}
    S(z) = I + O\left(\frac{1}{n}\right)
    \end{equation}
uniformly for $z \in \mathbb C \setminus \Gamma_S$ and uniformly
for $A$ and $B$ in compact subsets of the set {\rm(\ref{ABset})}.
\end{proposition}

The estimate (\ref{asympS}) is the basic asymptotic result. Unravelling
the sequence of transformations $Y \mapsto U \mapsto \tilde{T} \mapsto T \mapsto S$,
we obtain asymptotic formulas for $Y$ in any region of the
complex plane. In this way we obtain the asymptotic formulas for
$\widehat{P}_n = Y_{11}$.

\subsection{Proof of Theorem \ref{maintheorem}}

\begin{varproof} \textbf{of Theorem \ref{maintheorem}.}
Suppose $A$ and $B$ satisfy (\ref{ABinequalities}) and let $n \in \mathbb N$
such that $An$, $Bn$, and $(A+B)n$ are non-integers.
Then we have $\widehat{P}_n^{(An,Bn)}(z) = Y_{11}(z)$ by (\ref{formulaY}).

For $z$ in domain I away from the branch points, we get by using
 (\ref{defUz}), (\ref{tildeT1}), (\ref{defT}) and (\ref{defS1}),
\begin{eqnarray} \nonumber
    Y_{11}(z) &= & \left( \frac{G(z)}{\kappa} \right)^n U_{11}(z) \\
    & = & \nonumber
        \left( \frac{G(z)}{\kappa} \right)^n
        \left( \tilde{T}_{11}(z) - \frac{1}{d_L} e^{-2n\phi(z)} \tilde{T}_{12}(z) \right) \\
    & = & \nonumber
         \left( \frac{G(z)e^{-\phi(z)}}{\kappa} \right)^n
        \left( e^{n \phi(z)}T_{11}(z) - \frac{d_R}{d_C} e^{-n\phi(z)} T_{12}(z) \right) \\
    & = & \label{Y11formula1}
     \left( \frac{G(z)e^{-\phi(z)}}{\kappa} \right)^n
        \left( e^{n\phi(z)}(SN)_{11}(z) - \frac{d_R}{d_C} e^{-n\phi(z)} (SN)_{12}(z) \right).
\end{eqnarray}

Since $S = I + O(\frac{1}{n})$ and since the entries of $N$ are bounded and bounded away from
zero away from the branch points, we have that
\begin{equation} \label{SN11andSN12}
    (SN)_{11} = N_{11}(I + O(1/n)) \quad \text{ and } (SN)_{12} = N_{12}(I + O(1/n)).
\end{equation}

Next we recall that for $z$ in domain I,
\[ \frac{G'(z)}{G(z)} = g'(z) = \frac{A+B+2}{2} \frac{R(z)}{z^2-1} - \frac{A/2}{z-1} - \frac{B/2}{z+1}, \]
so that
\[ \log G(z) = \phi(z) - \frac{A}{2} \log(z-1) - \frac{B}{2} \log(z+1)  + const. \]
Thus there is a constant $c$ such that
\begin{equation} \label{relationc}
     \frac{G(z)e^{-\phi(z)}}{\kappa} = e^{-c} (z-1)^{-A/2} (z+1)^{-B/2}
     \qquad \text{for $z$ in domain I}.
\end{equation}
Since $Y_{11}(z)$ is a monic polynomial of degree $n$, we find by letting $z \to \infty$
in (\ref{Y11formula1}) and using (\ref{defc}) and (\ref{relationc}),
that $c$ should be as defined in (\ref{defc}).

Plugging (\ref{SN11andSN12}), (\ref{relationc}) and the formulas (\ref{constants})
for $d_R$ and $d_C$ into formula (\ref{Y11formula1}) we obtain
(\ref{asformleft}) for $z$ in domain I.

\medskip
For $z$ in domain II away from the branch points, we find in the same way
\begin{eqnarray}
    Y_{11}(z) \label{Y11formula2}
    & = & \left( \frac{G(z)e^{\phi(z)}}{\kappa} \right)^n
        \left( e^{n\phi(z)} (SN)_{11}(z) - \frac{d_R}{d_C} e^{-n\phi(z)}(SN)_{12}(z) \right)
\end{eqnarray}
Since $G_+ = G_- e^{-2\phi}$ on $\Gamma_L$, see (\ref{G+divideG-}), we have that
$G e^\phi$ is the analytic continuation of $G e^{-\phi}$ into domain II. So we have
by (\ref{relationc})
\begin{equation} \label{relationc2}
    \frac{G(z) e^{\phi}(z)}{\kappa} = e^{-c} (z-1)^{-A/2} (z+1)^{B/2}
    \qquad \text{for $z$ in domain II}.
\end{equation}
Then using (\ref{constants}), (\ref{SN11andSN12}), and (\ref{relationc2})
in (\ref{Y11formula2}), we obtain (\ref{asformleft}) for $z$ in domain II.

\medskip
For $z$ in domain III away from the branch points, we obtain
\begin{eqnarray} \label{Y11formula3}
    Y_{11}(z) &= & \left( \frac{G(z)e^{\phi(z)}}{\kappa} \right)^n
        \left( - \frac{d_L}{d_C} e^{n\phi(z)} (SN)_{11}(z) - \frac{d_R}{d_C} e^{-n\phi(z)} (SN)_{12}(z) \right).
\end{eqnarray}
Again we use (\ref{constants}), (\ref{SN11andSN12}), and (\ref{relationc2})
to obtain (\ref{asformmiddle1}) from (\ref{Y11formula3}) for $z$ in domain III.
\medskip

The proofs of the formulas (\ref{asformmiddle2}) and (\ref{asformright})
 for $z$ in the other domains IV, V, and VI are the same.
\medskip

We have derived the formulas
(\ref{asformleft})--(\ref{asformright}) under the assumption that
$n$ is such that $An$, $Bn$, and $(A+B)n$ are non-integers. Since
the formulas hold uniformly in $A$ and $B$ in compact subsets of
(\ref{ABset}) and $\widehat{P}_n^{(An,Bn)}$ depends continuously
on $A$ and $B$, they continue to hold if $An$ or $Bn$ is an
integer. However, we cannot allow $(A+B)n$ to be an integer, since
then there is a reduction in the degree of $P_n^{(An,Bn)}$ and we
cannot normalize the Jacobi polynomial to be monic.

This completes the proof of Theorem \ref{maintheorem}.
\end{varproof}

\subsection{Proof of Theorem \ref{asymnearbranch}}

\begin{varproof} \textbf{of Theorem \ref{asymnearbranch}.}
Suppose $A$ and $B$ satisfy (\ref{ABinequalities}) and let $n \in \mathbb N$
such that $An$, $Bn$, and $(A+B)n$ are non-integers.
Then we have $\widehat{P}_n^{(An,Bn)}(z) = Y_{11}(z)$ by (\ref{formulaY}).

Let $z \in U_{\delta}$ be in domain VI. Following the transformations
(\ref{defUz}), (\ref{tildeT6}), (\ref{defT}), we see that
\[ \begin{pmatrix} Y_{11}(z) \\ * \end{pmatrix} =
    \left(\frac{G(z) e^{-\phi(z)}}{\kappa} \right)^{n}
    T(z) \begin{pmatrix} e^{n\phi(z)} \\ - \frac{d_L}{d_C} e^{-n\phi(z)}
        \end{pmatrix},
\]
where $*$ denotes an unimportant unspecified entry.
For $z$ in domain VI, we have (\ref{relationc}) so that
\[\begin{pmatrix} Y_{11}(z) \\ * \end{pmatrix} =
    e^{-nc} (z-1)^{-An/2} (z+1)^{-Bn/2}
    T(z) \begin{pmatrix} e^{n\phi(z)} \\ - \frac{d_L}{d_C} e^{-n\phi(z)}
        \end{pmatrix}.
\]
Since $z$ belongs to $U_{\delta}$, we have $T(z) = S(z)P(z)$ by
(\ref{defS2}). By (\ref{defP}) we have  $P(z) =
E^{-t}(z) \Psi^{-t}(s) e^{-n\phi(z) \sigma_3}$ where $s = n^{2/3} f(z)$. Thus
\begin{equation} \label{Y11nearbranch}
    \begin{pmatrix} Y_{11}(z) \\ * \end{pmatrix} =
    e^{-nc} (z-1)^{-An/2} (z+1)^{-Bn/2}
    S(z)E^{-t}(z) \Psi^{-t}(s)\begin{pmatrix} 1 \\ - \frac{d_L}{d_C}
        \end{pmatrix}.
\end{equation}

From (\ref{defE}) we see that
\begin{equation} \label{defE2}
    E^{-t}(z)=\frac{1}{2\sqrt{\pi}}\, e^{-\pi i/6} \begin{pmatrix} 1 & -1
    \\ i & i \end{pmatrix}\, \left( \frac{a(z)}{s^{1/4}}\right)^{\sigma_3}.
\end{equation}
Furthermore, we have $0 < \arg s < \pi/3$ for $s = n^{2/3} f(z)$ since
$z$ is in domain VI, so that we  use the formula (\ref{defPsi1}) to
evaluate $\Psi^{-t}(s)$. Taking into account that
\[
\det \left(\begin{array}{cc}
    \Ai(s) & \Ai(\omega^2 s) \\
    \Ai'(s) & \omega^2 \Ai'(\omega^2 s) \end{array} \right) =\frac{1}{2\pi
    }\, e^{\pi i/6},
\]
we have from (\ref{defE2}) and (\ref{defPsi1}),
\[ E^{-t}(z) \Psi^{-t}(s) =
    \sqrt{\pi} e^{\pi i/6}
    \begin{pmatrix} -i & i \\ 1 & 1 \end{pmatrix}
    \left( \frac{a(z)}{s^{1/4}}\right)^{\sigma_3}
    \begin{pmatrix} \omega^2 \Ai'(\omega^2 s) & - \Ai'(s) \\
    - \Ai(\omega^2 s) & \Ai(s) \end{pmatrix} e^{\pi i/6 \sigma_3}. \]
Plugging this into (\ref{Y11nearbranch}) we get
\[ \begin{pmatrix} Y_{11}(z) \\ * \end{pmatrix} =
   e^{-nc} (z-1)^{-An/2}(z+1)^{-Bn/2}
    \sqrt{\pi} S(z) \begin{pmatrix} -i & i \\ 1 & 1 \end{pmatrix}
    \left( \frac{a(z)}{s^{1/4}} \right)^{\sigma_3}
    \begin{pmatrix} - \left(-\frac{d_L}{d_C} \Ai(s) + \omega^2 \Ai(\omega^2 s)\right)'    \\
    -\frac{d_L}{d_C} \Ai(s) + \omega^2 \Ai(\omega^2 s) \end{pmatrix}
\]
where the prime denotes the derivative with respect to $s$.

Comparing with (\ref{constants}) and (\ref{Airycombination}) we
see that
\[ \A(s) = \A(s;A,B,n) = - \frac{d_L}{d_C} \Ai(s) + \omega^2 \Ai(\omega^2s). \]
Thus
\begin{equation}\label{asympLocal}
\begin{pmatrix}
Y_{11}(z) \\ * \end{pmatrix} =
    e^{-nc} (z-1)^{-An/2} (z+1)^{-Bn/2}
    \sqrt{\pi} S(z) \begin{pmatrix} i \frac{s^{1/4}}{a(z)} \A(s) + i \frac{a(z)}{s^{1/4}} \A'(s) \\[10pt]
        \frac{s^{1/4}}{a(z)} \A(s) - \frac{a(z)}{s^{1/4}} \A'(s)
        \end{pmatrix}
\end{equation}
We derived the formula (\ref{asympLocal}) for $z \in U_{\delta}$ in the domain VI.
Similar calculations for $z \in U_{\delta}$ in
the other domains give the same result, so (\ref{asympLocal}) is
valid in the full neighborhood $U_{\delta}$ of $\zeta_-$.
Now it remains to recall that $Y_{11} = \widetilde{P}_n^{(An,Bn)}$
and that $S = I + O(1/n)$ as $n \to \infty$ in order to obtain
(\ref{asformnearzeta-}).

We have derived (\ref{asformnearzeta-}) under the assumption
that $An$ and $Bn$ are non-integers. Since the formula holds uniformly
for $A$ and $B$ in compact subsets of (\ref{ABset})
and $\widetilde{P}_n^{(An,Bn)}$ depends continuously on $A$
and $B$, they continue to hold if $An$ or $Bn$ is an integer.

This completes the proof of Theorem \ref{asymnearbranch}.
\end{varproof}

\subsection{Proof of Theorem \ref{theoremweak2}}

\begin{varproof} \textbf{of Theorem \ref{theoremweak2}.}
Conclusions of Theorem \ref{theoremweak2} are based upon the
strong asymptotics obtained in Theorem \ref{maintheorem}.

We let $(\alpha_n)$ and $(\beta_n)$ be two sequences such that
$\alpha/n \to A$ and $\beta_n/n \to B$ where $A$ and $B$ satisfy
the inequalities (\ref{ABinequalities}) and we assume that the
limits (\ref{closetoZ1})--(\ref{closetoZ3}) exist.

Taking into account that formula (\ref{asformleft}) is uniform in
$A$, $B$, for $z$ in domains I and II  we have
\begin{equation*}
\begin{split}
\widehat P_n^{(\alpha _n, \beta _n)}(z) = & \
e^{-nc} (z-1)^{-\alpha_n/2} (z+1)^{-\beta_n/2} \\
    &    \left(e^{n\phi(z)} N_{11}(z) \left(1+O\left(\frac{1}{n}\right)\right) -
        e^{-\alpha_n\pi i} \frac{\sin (\beta_n \pi)}{\sin((\alpha_n+\beta_n)
        \pi)}\,
        e^{-n\phi(z)} N_{12}(z)\left(1+O\left(\frac{1}{n}\right)\right)
        \right).
\end{split}
\end{equation*}
Since the first factors in the right hand side have no zeros
in domains I--II, $z$ is a zero of $P_n^{(\alpha _n, \beta _n)}$ only if
\begin{equation}\label{asymFor Zeros}
e^{2n\phi(z)}=
e^{-\alpha_n \pi i} \frac{\sin (\beta_n\pi)}{\sin((\alpha_n+\beta_n) \pi)}
 \frac{N_{12}(z)}{N_{11}(z)} \left(1+O\left(\frac{1}{n}\right)\right).
\end{equation}
But $N_{12}/N_{11}$ is uniformly bounded and uniformly bounded
away from zero, if we stay away from the branch points
$\zeta_\pm$. Thus, taking the absolute values in (\ref{asymFor
Zeros}), we see that the zeros in domains I--II away from the
branch points must satisfy
\[
2 \Re \phi(z) = \frac{1}{n} \log \left|
    \frac{\sin (\beta_n \pi)}{\sin((\alpha_n+\beta_n)\pi)}\right| +
    O\left(\frac{1}{n}\right).
\]

Analogously, the zeros of $P_n^{(\alpha_n, \beta_n)}$ in
the other domains III--VI away from
the branch points satisfy
\[ 2 \Re \phi(z) = \frac{1}{n}
    \log \left| \frac{ \sin(\beta_n \pi)}{ \sin(\alpha_n \pi)}\right| +
    O\left(\frac{1}{n} \right),
    \qquad \textrm{for $z$ in domain III}, \]
\[ 2 \Re \phi(z) = \frac{1}{n}
    \log \left| \frac{\sin(\alpha_n \pi)}{\sin(\beta_n \pi)} \right|
    + O\left(\frac{1}{n} \right),
    \qquad \textrm{for $z$ in domain IV}, \]
\[
  2 \Re \phi(z)=\frac{1}{n} \log \left|\frac{\sin(\alpha_n \pi)}{\sin((\alpha_n+\beta_n)\pi)}\right| +
    O\left(\frac{1}{n}\right),
    \qquad \textrm{for $z$ in domains V--VI}.
\]

Furthermore, for any sequence of real numbers $(\kappa_n)$,
$$
\lim_n \frac{1}{n}\, \log \left|\sin (\pi \kappa _n) \right|=
\lim_n \frac{1}{n}\, \log \left|\dist ( \kappa _n, \Z) \right|\,,
$$
whenever either one of these limits exists. Thus, under the
assumptions of the theorem, the zeros of $P_n^{(\alpha_n,
\beta_n)}$ away from the branch points must satisfy
\begin{equation}
\label{casesZeros} 2 \Re \phi(z) = r+
    O\left(\frac{1}{n}\right), \quad r=\begin{cases}  r_{\alpha
+\beta }-r_{\beta },
& \text{for $z$ in domains I and II,}\\
r_{\alpha  }-r_{\beta }, & \text{for $z$ in domain III,}\\
r_{\beta }-r_{\alpha },
& \text{for $z$ in domain IV,}\\
r_{\alpha +\beta }-r_{ \alpha  }, & \text{for $z$ in domains V and
VI.}
\end{cases}
\end{equation}
From the definition of the constants
(\ref{closetoZ1})--(\ref{closetoZ3}) it follows that the
``generic'' case is
\begin{equation}\label{genericCase}
r_{\alpha  }=r_{\beta }=r_{\alpha +\beta }\,.
\end{equation}
Recall that $\Re \phi(z)>0$ in domains I and VI, and $\Re
\phi(z)<0$ in domains II--V. Hence by (\ref{casesZeros}), if $z\in
\C\setminus \Gamma$, then $P_n^{(\alpha_n, \beta_n)}(z)\neq 0$ for
sufficiently large $n$, which proves that the zeros can accumulate
only on $\Gamma$.
%
%
%

Next assume we are in case (b) of the Theorem \ref{theoremweak2},
that is,
$$
r_{\alpha +\beta }>r_\alpha =r_\beta\,.
$$
By (\ref{casesZeros}), the zeros cannot accumulate in domains II,
III, IV and V, nor on $\Gamma_L\cup \Gamma_R$. In domains I and VI
they must satisfy
$$
 \Re \phi(z) = r+
    O\left(\frac{1}{n}\right), \quad r= \frac{ r_{\alpha
+\beta }-r_{\beta }}{2}=\frac{ r_{\alpha +\beta }-r_{ \alpha
}}{2}\,,
$$
showing that they must accumulate on the curve $\Gamma_r$ defined
in (\ref{defGammar}). Hence, in this case the accumulation set
belongs to $\Gamma_C \cup \Gamma_r$.

The rest of the cases is analyzed in the same fashion.

Once we have established that the zeros accumulate along curves in
the complex plane, it remains to find the asymptotic zero
distribution. We can use  the differential equation (see e.g.\
\cite[\S 4.22]{szego:1975})
$$
(1-z^2)\, y_n''(z) +\left[ \beta _n-\alpha _n-(\alpha _n+\beta
_n+2) z \right]\, y_n'(z) + n (n+\alpha _n+\beta _n+1)\, y_n(z)=0
$$
satisfied by $y_n=P_n^{(\alpha_n, \beta_n)}$. Rewriting this
equation in terms of $h_n=y'_n/(n y_n)$ we reduce it to the
Riccati form
\begin{equation}\label{riccati}
(1-z^2)\, \left( \frac{1}{n}\, h_n'(z)+h_n^2(z) \right) +
\frac{\beta _n-\alpha _n-(\alpha _n+\beta _n+2) z }{n} \, h_n'(z)
+ \frac{ n+\alpha _n+\beta _n+1 }{n}=0\,.
\end{equation}
Let $\nu_n$ denote the normalized zero counting measures of
$P_n^{(\alpha _n, \beta _n)}$. Using the week compactness of the
sequence $( \nu_n )$ we may take a subsequence $\Lambda\subset \N$
such that $\nu_n$ converge along $\Lambda$ to a certain unit
measure $\nu$ in the weak*-topology. By the discussion above,
$\nu$ is supported on a finite union of analytic arcs or curves
(level sets $\Gamma_r$), and every compact subset of $\C\setminus
\supp (\nu)$ contains no zeros of $P_n^{(\alpha _n, \beta _n)}$
for $n$ sufficiently large. Hence,
$$
h_n(z)=\int_\Gamma \frac{d\nu_n(t)}{z-t} \longrightarrow
h(z)=\int_\Gamma \frac{d\nu(t)}{z-t}\,, \quad n\in \Lambda\,,
$$
locally uniformly in $\C\setminus \supp (\nu)$. Taking limits in
(\ref{riccati}) we obtain that $h$ satisfies the quadratic
equation
$$
(1-z^2)\, h^2(z)  + \left[ B-A-(A+B) z \right]\, h(z) + A +B+1
=0\,,
$$
so that
$$
\int_\Gamma \frac{d\nu(t)}{z-t}=\frac{A+B+2}{2}\,
\frac{R(z)}{z^2-1}-\frac{1}{2}\, \left( \frac{A}{z-1}+
\frac{B}{z+1}\right)\,, \quad z \in \C\setminus \supp (\nu)\,,
$$
with $R$ defined in (\ref{defRz}) and $\zeta_\pm$ given in
(\ref{zetapm}). By Sokhotsky-Plemelj's formulas, on every arc of
$\supp (\nu)$,
\begin{equation}\label{nu}
   d\nu(z) = \frac{A + B+2}{2 \pi i} \frac{R_+(z)}{z^2-1} \,dz\,.
\end{equation}

Assume that (\ref{genericCase}) holds, so that $\supp (\nu)\subset
\Gamma$. By (\ref{nu}), $\nu'=\mu'$ almost everywhere on $\supp
(\nu)$. Taking into account Lemma \ref{lemma2} and that $\nu$ is a
unit measure it follows that $\nu=\mu$.

If $r_{\alpha +\beta }>r_\alpha =r_\beta$ then $\supp (\nu)\subset
\Gamma_C \cup \Gamma_r$, $r=(r_{\alpha +\beta }-r_{ \alpha
})/2>0$. Again taking into account Lemma \ref{lemma2} and the
normalization of $\nu$ it follows that $\supp (\nu)= \Gamma_C \cup
\Gamma_r$.

The remaining cases are analyzed analogously. This completes the
proof of Theorem \ref{theoremweak2}.
\end{varproof}

\section*{Acknowledgements}

This research was partially supported by the Ministry of Science
and Technology of Spain through the grant BFM2001-3878-C02-02.
A.B.J.K.\ was also supported by FWO research projects G.0176.02
and G.0455.04. Additionally, A.M.F.\ acknowledges the support of
Junta de Andaluc{\'\i}a, Grupo de Investigaci{\'o}n FQM 0229, and
of the Ministry of Education, Culture and Sports of Spain through
the grant PR2003--0104.

\end{document}